 \documentstyle[epsf,amssymb,12pt]{amsart}
 \input xy
 \xyoption{all}

 \newcommand{\la}{\langle}
 \newcommand{\ra}{\rangle}

 \newtheorem{theorem}{\bf Theorem}[section]
 \newtheorem{lemma}[theorem]{\bf Lemma}
 
 \newtheorem{remark}[theorem]{\bf Remark}
 \newtheorem{prop}[theorem]{\bf Proposition}
 \newtheorem{corollary}[theorem]{\bf Corollary}

 \DeclareFontFamily{OT1}{rsfs}{}
 \DeclareFontShape{OT1}{rsfs}{n}{it}{<->rsfs10}{}
 \DeclareMathAlphabet{\curly}{OT1}{rsfs}{n}{it}

 \newcommand{\CC}{{\Bbb C}}
 \newcommand{\DD}{{\Bbb D}}
 \newcommand{\EE}{{\Bbb E}}
 \newcommand{\HH}{{\Bbb H}}

 \newcommand{\PP}{{\Bbb P}}
 
 \newcommand{\RR}{{\Bbb R}}
 \newcommand{\VV}{{\Bbb V}}
 
 \newcommand{\ZZ}{{\Bbb Z}}
 \newcommand{\bV}{{\bf V}}

 \newcommand{\ulie}{{\frak u}}

\newcommand{\splie}{{\frak sp}}

 \newcommand{\alg}{\operatorname{alg}}
 \newcommand{\Aut}{\operatorname{Aut}}

 \newcommand{\Coh}{\operatorname{Coh}}
 
 \newcommand{\diag}{\operatorname{diag}}

 \newcommand{\End}{\operatorname{End}}

 \newcommand{\Ext}{\operatorname{Ext}}
 \newcommand{\free}{\operatorname{free}}
 \newcommand{\gen}{\operatorname{gen}}
 \newcommand{\GL}{\operatorname{GL}}

 \newcommand{\Hom}{\operatorname{Hom}}
 \newcommand{\Id}{\operatorname{Id}}
 
 \newcommand{\Jac}{\operatorname{Jac}}
 \newcommand{\Ker}{\operatorname{Ker}}
 
 \newcommand{\Lie}{\operatorname{Lie}}

 \newcommand{\mult}{\operatorname{mult}}

 \newcommand{\Pic}{\operatorname{Pic}}

 \newcommand{\rk}{\operatorname{rk}}

 \renewcommand{\Sp}{\operatorname{Sp}}

 \newcommand{\U}{\operatorname{U}}

 \renewcommand{\Im}{\operatorname{Im}}
 
 \newcommand{\ov}{\overline}
 \newcommand{\un}{\underline}

 \newcommand{\AAA}{{\curly A}}
 \newcommand{\BBB}{{\curly B}}
 \newcommand{\CCC}{{\curly C}}
 
 \newcommand{\EEE}{{\curly E}}
 \newcommand{\FFF}{{\curly F}}
 \newcommand{\GGG}{{\curly G}}
 \newcommand{\HHH}{{\curly H}}
 \newcommand{\III}{{\cal I}}
 
 \newcommand{\LLL}{{\curly L}}
 \newcommand{\MMM}{{\cal M}}
 \newcommand{\NNN}{{\cal N}}
 \newcommand{\cO}{{\cal O}}
 \newcommand{\OOO}{{\curly O}}
 \newcommand{\PPP}{{\curly P}}
 \newcommand{\QQQ}{{\curly Q}}
 \newcommand{\RRR}{{\cal R}}
 \newcommand{\SSS}{{\cal S}}
 \newcommand{\TTT}{{\curly T}}
 
 \newcommand{\VVV}{{\curly V}}
 \newcommand{\XXX}{{\cal X}}

 \newcommand{\qu}{/\kern-.7ex/}
 \newcommand{\exh}{\to\kern-1.8ex\to}

 \newcommand{\uC}{\un{\CC}}

 \newcommand{\VP}{{\curly V}\kern-0.9ex\PPP}

 \newcommand{\imag}{{\mathbf i}}

 \newcommand{\bx}{{\mathbf x}}

\newcommand{\Beta}{B}

 \title[Representations in $\Sp(4,\RR)$]
 {Representations of the fundamental group
 of a closed oriented surface in $\Sp(4,\RR)$}

\author{O. Garc\'{\i}a--Prada, I. Mundet i Riera}

\date{21 March  2003}

\begin{document}
 \maketitle


 \section{Introduction}

Let $X$ be a closed oriented surface of genus $g\geq 2$ 
and let $\pi=\pi_1(X)$ be  its
fundamental group.  Let $\Sp(4,\RR)$ be  the group of linear 
transformations of $\RR^4$ preserving its standard real symplectic form. 
Consider the set $\XXX:=\Hom(\pi,\Sp(4,\RR))$  of group homomorphisms
from  $\pi$ to $\Sp(4,\RR)$. We also refer to the elements of $\XXX$ 
as {\bf representations} of $\pi$ in $\Sp(4,\RR)$.
 Using a set of  generators of $\pi$, 
$\XXX$ can be embedded in $\Sp(4,\RR)^{2g}$,  acquiring in this
way a natural   topological structure.
The goal of this paper is to compute the number of connected components
of $\XXX$.

Given a representation of $\pi$ in $\Sp(4,\RR)$, there is an integer  
associated to it, which is geometrically obtained by considering
the flat $\Sp(4,\RR)$-bundle corresponding to the representation and taking
a reduction of the structure group to $\U(2)$ --- the maximal compact 
subgroup of $\Sp(4,\RR)$. The {\bf degree} $d$ of this $\U(2)$-bundle is an invariant 
of the  representation, which we call the degree of the representation.
 Let $\XXX(d)\subset \XXX$ be the set of 
representations of degree $d$.
There is a Milnor-Wood type inequality   (\cite{DT,Tu}) which says that 
$\XXX(d)$ is empty unless  $|d|\leq 2g-2$.  Our  main result in  this paper
is the following.

\begin{theorem} \label{main-theorem}
Let $d$ be an integer such that  $|d|\leq 2g-2$. Then

(1) $\XXX(d)$ is non-empty and connected if $|d|< 2g-2$;

(2) $\XXX(d)$  has $3.2^{2g}+2g-4$ non-empty connected components
if  $|d|= 2g-2$.
\end{theorem} 

The proof of (2) and the case $d=0$ for the space of 
 reductive representations is due to Gothen \cite{Go2}. 
Recall that a representation $\rho$  from $\pi$ to a real algebraic
group $G$  is called {\bf reductive} if the real Zariski closure of $\rho(\pi)$ 
in $G$ is a reductive group. When the closure coincides with $G$, the 
representation
is said to be {\bf irreducible}. 
Let $\XXX^+\subset \XXX$ be the subspace  of reductive representations 
of $\pi$ in  $\Sp(4,\RR)$.  In Section \ref{sec-rep}  we show that the 
inclusion $\XXX^+\subset \XXX$ induces a bijection of connected
components (Theorem \ref{red-versus-nonred}). In fact,  we show, this  result is 
valid replacing  $\Sp(4,\RR)$ by any reductive real algebraic group $G$ and  
$\pi$ by any  finitely generated group.

Let $\XXX^+(d)$ be the space  of reductive representations from $\pi$ to
$\Sp(4,\RR)$ of degree $d$. This space is invariant under the  action of $\Sp(4,\RR)$ 
by conjugation. The {\bf moduli space of representations} of degree $d$ 
is defined as $\RRR(d):=\XXX^+(d)/\Sp(4,\RR)$. The reductivity condition is
precisely what is needed in order for  $\RRR(d)$ to be a  Hausdorff space.
 Gothen's results can now be stated as follows.
\begin{theorem}[Gothen \cite{Go2}]\label{gothen-thm}
The moduli space  $\RRR(0)$ is non-empty and connected, and  $\RRR(\pm(2g-2))$ has $3.2^{2g}+ 2g-4$
non-empty connected components.
\end{theorem}
In this paper we  settle  the situation for the
remaining cases.  More precisely, we prove the following.  
\begin{theorem}
\label{thm:mainthm-red}
Let $d$ be any integer satisfying $0<|d|<2g-2$.
The moduli space $\RRR(d)$ is connected and the subspace of irreducible 
representations is non-empty.
\end{theorem}

To prove this theorem, in Section \ref{sec-higgs}, 
we follow the methods of Gothen (\cite{Go2}) and
Hitchin (\cite{Hi1,Hi2}) by choosing a complex structure on $X$ and exploting the relation 
between   $\RRR(d)$ and the {\bf moduli space of polystable
$\Sp(4,\RR)$-Higgs bundles} $\MMM(d)$.
 These are Higgs bundles $(E,\Phi)$, 
where  
$E=V\oplus V^*$, for $V$ a rank 2 holomorphic vector
bundle of degree $d$, and  Higgs field $\Phi: E\to E^\ast\otimes K$
 of the form
\begin{equation}
\Phi=\left(
\begin{array}{cc}
0 & \beta \\ \gamma & 0\end{array}\right)
   :V\oplus V^*\to  (V\oplus V^*)\otimes K
\end{equation}
with $\beta\in H^0(S^2 V\otimes K)$ and $\gamma \in H^0(S^2 V^*\otimes K)$.
The basic fact
is that $\MMM(d)$ is a 
complex analytic variety which is homeomorphic to $\RRR(d)$ (by results of
Hitchin, Donaldson, Simpson and  Corlette).
This homeomorphism induces a homeomorphism between the subspace of
irreducible  representation in $\RRR(d)$ and the subspace
of stable Higgs bundles in $\MMM(d)$.

By solving Hitchin equations
for a Hermitian metric on $V$ and considering the square of the $L^2$-norm 
of $\Phi$, one has  a proper function on $\MMM(d)$.
The number of connected components of the local minima of this proper
function gives an upper bound on
the number of connected components of $\MMM(d)$ and hence of 
$\RRR(d)$.
These local
minima have been characterised by Gothen \cite{Go2} for any $d$. 
In particular when $|d|<2g-2$, Gothen shows that the minima coincide
with the subvariety $\NNN(d)\subset \MMM(d)$ for which either
$\beta=0$ or $\gamma=0$ (Proposition \ref{minima}). Which one of the sections
 actually vanishes is determined by
the sign of $d$. We prove that if $0<|d|<2g-2$, $\NNN(d)$ is connected.
To show this,  in Section \ref{sec-pairs} we study  a  more general 
situation that is  of independent interest. Namely, we introduce  a 
(poly)stability
criterium for pairs $(V,\beta)$, where $V$ is a rank 2 holomorphic 
vector bundle of degree $d$, and 
$\beta\in H^0(S^2 V\otimes K)$ which  depends on a real parameter
$\alpha$,   thus defining moduli spaces $\NNN_\alpha(d)$. It turns out 
that the subvariety of minima $\NNN(d)$ can be identified with
$\NNN_0(d)$ (Proposition \ref{prop:stabd}). Our main result 
concerning $\NNN_\alpha(d)$ is 
the following Theorem proved in Section \ref{sec-conn}.

\begin{theorem}\label{thm:conn}
Let $-2(g-1)< d < 0$ be any integer and let $\alpha\geq 0$ be a real
number. The moduli space $\NNN_\alpha(d)$ is connected. Moreover, the 
subvariety consisting of stable pairs is non-empty.
\end{theorem}

From this we deduce the following.
\begin{theorem}\label{higgs-moduli}
Let $d$ be any integer such that $0<|d|< 2g-2$. The moduli space
$\MMM(d)$ of polystable $\Sp(4,\RR)$-Higgs bundles of degree $d$ is connected. 
Moreover, the  subvariety consisting of stable Higgs bundles is non-empty.
\end{theorem}

This proves Theorem \ref{thm:mainthm-red}, which combined  with Gothen's 
Theorem  
\ref{gothen-thm} and  Theorem \ref{red-versus-nonred} proves  Theorem 
\ref{main-theorem},
since the connectedness of $\Sp(4,\RR)$ implies that the number of connected 
components of $\XXX^+(d)$ coincides with that of $\RRR(d)$.

\section{Representations of the fundamental group}\label{sec-rep}

\subsection{Reductive representations of $\pi_1(X)$}

The results of this section and the next one  apply also when
$\pi$ is replaced by any finitely generated group.

Let $X$ be a closed oriented surface of genus $g\geq 2$ 
and let $\pi=\pi_1(X)$ denote its
fundamental group.
Let $G\subset \GL(N,\RR)$ be a non-compact  real reductive algebraic 
group (when $G$ is compact everything which follows in this section
and the next one are also true for obvious reasons). Let  
$\XXX:=\Hom(\pi,G)$ be the set of representations of $\pi$
 in $G$, and let 
 $$\XXX^+:=\{\rho\in\XXX\mid \ov{\rho(\pi)}\subset G
 \text{ is reductive }\}$$
 the set of reductive representations ($\ov{\rho(\pi)}$
 denotes the real Zariski closure of $\rho(\pi)$ in $G$).
 Take generators 
 $\gamma_1,\dots,\gamma_k$ of $\pi$ and consider the inclusion
 $j:\XXX\to \EE:=\End(\RR^N)^k$ which sends $\rho$ to 
 $(\rho(\gamma_1),\dots,\rho(\gamma_k))\in G^k\subset\End(\RR^N)^k$.
 Consider  on $\XXX$ the topology induced by $j$ and the standard
 topology on the vector space $\EE$.
This is independent of the choice of generators of $\pi$. We take  
on $\XXX^+$ the  topology induced by the inclusion $\XXX^+\subset\XXX$.
Furthermore,
 $j(\XXX)\subset \EE$ is a real algebraic 
(affine) subvariety of
$G^k\subset\End(\RR^N)^k$, whose equations are defined by  requiring
the coordinates of a point $(\rho_1,\dots,\rho_k)\in G^k$ to
satisfy any relation satisfied by the generators $\{\gamma_j\}$. 
It follows by a theorem of Whitney \cite{Wh} that $j(\XXX)$ (and hence $\XXX$)
has a finite number of connected components.

 Consider the adjoint action of $G$ on $\XXX$: if $g\in G$ and
 $\rho\in\XXX$ then $g\cdot\rho$ is the representation defined by
 $g\cdot\rho(\gamma):=g \rho(\gamma) g^{-1}$ for any $\gamma\in\pi$.
 We consider similarly the diagonal adjoint action of $G$ on
 $\EE=(\End \RR^N)^k$, in such a way that the inclusion $j$ is $G$-equivariant.
 One has the following.
 
 \begin{theorem}[Richardson, Theorem 11.4 in \cite{Ri}]
 \label{thm:richardson}
 A representation $\rho\in\XXX$ is reductive if and only if
 the orbit $G\cdot j(\rho)\subset \EE$ is closed in the usual topology
 of $\EE$.
 \end{theorem}

 \begin{remark} 
 If a real algebraic group acts linearly on a vector space
 then an orbit which is closed in the usual topology may fail to be
 closed in the real Zariski topology (in general its Zariski closure
 will consist of a finite number of orbits which are closed in the
 usual topology). 
 This is in contrast to the situation over the complex numbers, where
 an orbit of a linear action is closed in the usual topology if and
 only if it is closed in the Zariski topology.
 \end{remark}

It follows from Theorem \ref{thm:richardson} that $\XXX^+/G$ with the
quotient topology is a Hausdorff space. 
The space $\RRR:=\XXX^+/G$ is called  the 
{\bf moduli space of representations} of $\pi$ in $G$.

Our main concern is the  study of  the connectedness of 
$\XXX$, $\XXX^+$ and $\RRR$. Of course, when $G$ is connected the number
of connected components of $\RRR$ coincides with that of $\XXX^+$.

 \subsection{From reductive representations to arbitrary representations}

 \begin{theorem}\label{red-versus-nonred}
 The inclusion $i:\XXX^+\subset\XXX$
 induces a bijection of connected components.
 \end{theorem}
 \begin{pf}

To prove this we need the following.
(All topological notions refer to the usual topology in $\EE$.)

 \begin{theorem}[Luna, Theorem 2.7 in \cite{Lu}]
 \label{thm:luna}
 Let $G$ be a real algebraic group acting linearly on a real vector space 
$\EE$.
 Let $\EE^+:=\{x\in \EE\mid G\cdot x\subset \EE\text{ is closed }\}$.
 For any $x\in \EE$ there is a unique closed orbit $p(x)$
 contained in the closure of $G\cdot x$. The space of orbits $\EE^+/G$ 
endowed with the quotient topology
 is Hausdorff, and the map  $p:\EE\to \EE^+/ G$ is continuous.
 \end{theorem}

 In \cite{Lu} this theorem is stated in a slightly different form.
 Luna defines a $G$-invariant subset $A\subset \EE$ to be $G$-saturated
 if it contains any orbit whose closure intersects the closure
 of any orbit inside $A$. Then he proves that for any $G$-invariant
 closed subset $F\subset \EE$ the smallest $G$-saturated subset of $\EE$
 which contains $F$ is also closed (this is Property (C) of \cite{Lu}).
 To deduce from this the continuity of $p$ note that if $C\subset
 \EE^+/G$ is closed then $p^{-1}(C)$ is the smallest $G$-saturated
 set containing the closure of $q^{-1}(C)\subset \EE$, where 
$q:\EE^+\to \EE^+/G$
 is the quotient map.

 Now we prove Theorem \ref{red-versus-nonred}. 
Denote by $G_0\subset G$ the connected 
 component of the identity. We first check that the map
 $i_*:\pi_0(\XXX^+)\to\pi_0(\XXX)$ is onto. For that, take any
 $\rho\in\XXX$. By Theorem \ref{thm:luna} the closure of
 $G_0\cdot j(\rho)$ contains a closed orbit $G_0\cdot j(\rho_0)$. On the
 other hand, since $G_0$ is connected, the closure of 
 $G_0\cdot j(\rho)$ is locally arc-connected, and $j$ is injective,
 it follows that there is a continuous map $c:[0,1]\to\XXX$ such that
 $c(0)=\rho_0$ and $c(1)=\rho$. Finally, Theorem \ref{thm:richardson}
 implies that $\rho_0\in\XXX^+$.

 To prove that $i_*$ is injective consider the following diagram,
 which commutes by the definition of $p$ (see Theorem \ref{thm:luna})
 $$\xymatrix{\XXX^+ \ar[r]^i \ar[rd]_q & \XXX \ar[d]^{p} \\
 & \XXX^+/G_0.}$$ 
 Since $G_0$ is connected it follows that $q_*:\pi_0(\XXX^+)\to
 \pi_0(\XXX^+/G_0)$ is a bijection, and since $q_*=p_*\circ i_*$
 we deduce that $i_*$ is an injection.
 \end{pf}

\subsection{Representations of $\pi_1(X)$ in $\Sp(4,\RR)$}
We use the same notation as in the previous two sections.
In particular Let $\EE=\End (\RR^4)^k$.
Let now $G=\Sp(4,\RR)$ and let  $\XXX:=\Hom(\pi,\Sp(4,\RR))$ be the set of 
representations of $\pi$ in $\Sp(4)$. 
To understand the geometric meaning of reductive representations of $\pi$ in 
$\Sp(4,\RR)$, consider the standard symplectic structure $\omega$ on $\RR^4$.
Recall that a subspace $\VV\subset\RR^4$ is called 
{\bf coisotropic} if 
$$
\VV^{\perp}=\{v\in\RR^4\mid
\omega(v,w)=0\ \forall w\in \VV\}\subset \VV.
$$
One has the following.

\begin{prop}
\label{prop:equivalencies}
Let $\rho:\pi\to\Sp(4,\RR)$ be a representation. 
The following are 
equivalent:

(i) $\rho$ is reductive;

(ii) $\Sp(4,\RR)\cdot j(\rho)\subset \EE$ is closed;

(iii) for any coisotropic
subspace $\VV\subset\RR^4$ which is preserved by $\rho$
there is a splitting 
\begin{equation}
\RR^4=\VV/\VV^{\perp}\oplus (\RR^4/\VV\oplus \VV^{\perp})
\label{eq:splitting}
\end{equation}
such that the image of $\rho$ is contained in 
$\Sp(\VV/\VV^{\perp})\times\Sp(\RR^4/\VV\oplus \VV^{\perp})$
(one can check that both $\VV/\VV^{\perp}$ and 
$\RR^4/\VV\oplus \VV^{\perp}$
carry natural symplectic structures).

\end{prop}

\begin{pf}
The equivalence between
(i) and (ii) is given by  Theorem \ref{thm:richardson}.
We just sketch the proof of the equivalence between (ii) and (iii) 
(the argument is similar to that of Lemma 28 in
\cite{Gu}).
Take a representation $\rho$ and suppose that there is
a sequence $\{g_j\}\subset\Sp(4,\RR)$ which is not
contained in any compact subset of $\Sp(4,\RR)$ and such that
for any $k$
$$\{g_j\rho(\gamma_k)g_j^{-1}\}\to\gamma_k'\in\Sp(4,\RR)
\qquad\text{as }j\to\infty.$$
Define $s_j:=g_j^*g_j$, where we take the adjoint of $g_j$
with respect to the standard Euclidean metric in $\RR^4$. Then
$s_j$ diagonialises in some basis $e_1,\dots,e_4$ with eigenvalues
$\lambda_1\geq\dots\geq\lambda_4$. Define
$\VV^j:=\la e_2,e_3,e_4\ra$ if $\lambda_1>\lambda_2$ or
$\VV^j:=\la e_3,e_4\ra$ otherwise.
Since $s_j\in\Sp(4,\RR)$, $\VV^j$ is coisotropic. Let
$\VV$ be a limit of a partial sequence of $\{\VV^j\}$. Then, since
$\{g_j\rho(\gamma_k)g_j^{-1}\}$ is bounded, $\rho$ preserves
$\VV$. Now, if $\gamma_k'$ belongs to the adjoint
orbit of $\gamma_k$ then one must have a splitting as in
(\ref{eq:splitting}). Finally, if the limit representation
defined by $\{\gamma_k'\}$ does not belong to the orbit of
$\rho$, then one can take a sequence $\{g_j\}$ in such a way
that no splitting of the form (\ref{eq:splitting})
is preserved by $\rho(\gamma_k)$ for all $k$.
\end{pf}

\begin{remark}
It is clear that a reductive representation $\rho\in\XXX$
is  {\bf irreducible}  if there is no nontrivial splitting
$\RR^4=\VV_1\oplus \VV_2$ with both $\VV_1$ and $\VV_2$ symplectic
subspaces of $\RR^4$ such that the image of $\rho$ is contained
in $\Sp(\VV_1)\times\Sp(\VV_2)$. 
\end{remark}

\subsection{Characteristic numbers} 

Our first step in determining  the number of connected 
components of $\XXX:= \Hom(\pi, \Sp(4,\RR))$
is  to study the topological invariants
associated to a representation.
For any $G$ there is a locally constant obstruction map
$$
o:\Hom(\pi,G)\to H^2(X,\pi_1(G)).
$$
For $G=\Sp(4,\RR)$, we have $\pi_1(G)\cong\ZZ$;
picking the standard isomorphism, the obstruction map
associates to any $\rho\in\XXX$ an integer $o(\rho)$ which we call 
the {\bf degree} of $\rho$, and which
can be obtained more geometrically by considering
the flat $\Sp(4,\RR)$ bundle associated to $\rho$.
This bundle admits
a reduction of the structure group to $\U(2)$ --- the maximal compact
subgroup of $\Sp(4,\RR)$ ---, whose first Chern class is precisely the degree  
of the representation. (Note  that the reduction
will not be in general compatible with the flat structure).
One  can also show  that this degree measures the  obstruction to the 
existence of a Lagrangian subbundle of the flat $\Sp(4,\RR)$-bundle.

Let $d$ be an integer; recall that
$\XXX(d)\subset \XXX$ consists of those 
representations whose degree is $d$. The subset $\XXX^+(d)\subset \XXX(d)$
is given by the reductive representations of  degree $d$ and 
$\RRR(d)\subset \RRR$ given by $\XXX^+(d)/\Sp(4,\RR)$ is the 
{\bf moduli space of representations of degree $d$}.
Our main goal is the study of the non-emptiness and connectedness of 
$\RRR(d)$. 

A restriction on the possible degrees of an element in $\XXX$ is 
given by the following Milnor-Wood type inequality(
Turaev \cite{Tu}, Domic and Toledo \cite{DT}).
\begin{prop}
Let $\rho$ be a representation of $\pi$ in $\Sp(4,\RR)$, and let $d$ be
the degree of $\rho$. Then
$$
|d|\leq 2g-2.
$$
\end{prop}
Hence $\XXX(d)$ is empty for $|d|> 2g-2$.

\section{$\Sp(4,\RR)$-representations of $\pi$ and Higgs bundles}
\label{sec-higgs}
In this section we follow  the methods of Hitchin \cite{Hi1,Hi2} and Gothen
\cite{Go1,Go2}. We refer to their papers for details.

\subsection{$\Sp(4,\RR)$-Higgs bundles}

\label{ss:sympHiggs}
As in the case of complex representations of $\pi$, representations
in $\Sp(4,\RR)$ are related to Higgs bundles.
To recall the basic ingredients  of this theory, we  fix from now on a
complex  structure on $X$.
A {\bf Higgs bundle} on $X$ is a pair
$(E,\Phi)$, where $E$ is a holomorphic vector bundle over $X$
and $\Phi: E\to E\otimes K$ is a holomorphic map, where $K$ is the canonical
line bundle of $X$.
The Higgs bundle $(E,\Phi)$ is said to be {\bf stable} if for any
proper subbundle $F\subset E$ such that $\Phi(F)\subset F\otimes K$
we have $\mu(F)<\mu(E)$, where $\mu(F)=\deg F/\rk F$.
The Higgs bundle $(E,\Phi)$ is said to be {\bf polystable}
 if it is the direct sum
of stable Higgs bundles of the same slope $\mu(E)$.
This condition appears as the requirement to solve Hitchin's equations.
More precisely, we have the following.
\begin{theorem}\cite{Hi1,Si1}
Let $(E,\Phi)$ be  a Higgs bundle with $\deg E=0$.
Then $E$ admits a Hermitian metric $H$  satisfying
\begin{equation}
F_H+[\Phi,\Phi^{*_H}]=0
\label{eq:EHh}
\end{equation}
if and only if $(E,\Phi)$ is polystable. (Here $F_H$ is the curvature
of the unique connection compatible with $H$ and the holomorphic 
structure of $E$.) Furthermore, the set of metrics $H$ which solve
(\ref{eq:EHh}) is convex, i.e., if $H,H'$ are both solutions then
for any real number $t\in [0,1]$ the metric $tH+(1-t)H'$ is also
a solution.
\end{theorem}

The last statement of the theorem follows from the fact that
if $(E,\Phi)$ is stable then there is a unique Hermitian metric $H$
solving (\ref{eq:EHh}).

The  particular class of Higgs bundles that will be of relevance in relation
to representations of $\pi$ in $\Sp(4,\RR)$ is given by pairs $(E,\Phi)$,
where  $E= V\oplus V^*$, with $V$ a rank 2 holomorphic vector bundle and
\begin{equation}
\Phi=\left(
\begin{array}{cc}
0 & \beta \\ \gamma & 0\end{array}\right)
   :V\oplus V^*\to (V\oplus V^*)\otimes K,
\label{eq:defPhi}
\end{equation}
where $\beta\in H^0(S^2 V\otimes K)$ and $\gamma\in H^0(S^2 V^*\otimes K)$.

A $\Sp(4,\RR)$-Higgs bundle can be regarded as a pair $(V,\varphi)$ 
where $V$ is a rank 2 holomorphic bundle and 
$\varphi=(\beta,\gamma)\in H^0( S^2 V\otimes K) \oplus H^0(S^2 V^*\otimes K)$.
Two $\Sp(4,\RR)$-Higgs bundles $(V,\varphi)$ and $(V',\varphi')$ are
{\bf isomorphic} if there is an isomorphism $\psi:V'\to V$ such that
$\varphi'=\psi^*\varphi$.
Let $d$ be an integer. Let $\MMM(d)$ be the moduli space of 
polystable $\Sp(4,\RR)$-Higgs bundles  $(V,\varphi)$ such that $\deg V=d$.
By stability of $(V,\varphi)$ we mean  stability of the corresponding
Higgs bundle $(E,\Phi)$, with $E= V\oplus V^*$, and $\Phi$ given by 
(\ref{eq:defPhi}).

The following refinement of Theorem \ref{eq:EHh} is necessary to relate 
the moduli space of $\Sp(4,\RR)$-Higgs bundles $\MMM(d)$ to $\RRR(d)$.

\begin{theorem}
\label{thm:HK}
A $\Sp(4,\RR)$-Higgs bundle $(V,\varphi)$ is polystable if and only if
there exists a Hermitian metric $h$ on $V$ such that
\begin{equation}
F_h+(\beta\beta^{*_h}-\gamma^{*_h}\gamma)=0.
\label{eq:HE0}
\end{equation}
where $F_h$ denotes the curvature of the unique connection 
compatible  with $h$ and the holomorphic structure of $V$.
\end{theorem}
\begin{pf}

Suppose that we have a metric
$h$ on $V$ for which (\ref{eq:HE0}) holds. Then, setting $H=h+h^*$
we get a solution to (\ref{eq:EHh}), so that $(E,\Phi)$ is polystable
and hence (by definition) $(V,\varphi)$ is also polystable.

To prove the converse, assume that $(V,\varphi)$ is polystable, i.e.
$(E,\Phi)$ is polystable. By Theorem \ref{eq:EHh} there is a metric $H$
on $E$ which solves (\ref{eq:EHh}).
We want to show  that this metric
can be taken of the form $H=h+h^*$ with respect to the splitting
$E=V\oplus V^*$.
We can view $H$ as a section of $E\otimes\ov{E}^*$, and using
the splitting 
$$E\otimes\ov{E}^*=
V\otimes \ov{V}^*\oplus V\otimes (\ov{V^*})^*
\oplus V^*\otimes \ov{V}^*\oplus V^*\otimes (\ov{V^*})^*$$
we may write $H=H_{00}+H_{01}+H_{10}+H_{11}$.
We claim that if $H':=H_{00}-H_{01}-H_{10}+H_{11}$ is also
a solution to (\ref{eq:EHh}). To prove the claim, consider
a local holomorphic framing $e_1,\dots,e_{2n}$ of $E$ whose
first (resp. last) $n$ sections give a framing of $V$ (resp. $V^*$).
Let $M_H$ be the $2n\times 2n$
matrix whose $(i,j)$ entry is $\la e_i,e_j\ra_H$.
Define $\rho:=\diag(\imag\Id_V,-\imag\Id_{V^*})$.
A simple computation shows that $M_{H'}=\rho M_H \rho^{-1}$.
On the other hand, if we identify $\Phi$ with a matrix by
means of our framing, then we have
$\Phi^{*_H}=(M_H^*)^{-1}\Phi^*M_H^*,$
from which we easily compute
\begin{equation}
[\Phi,\Phi^{*_{H'}}]=\rho [\rho^{-1}\Phi\rho,
(\rho^{-1}\Phi\rho)^{*_H}]\rho^{-1}=\rho[\Phi,\Phi^{*_H}]\rho^{-1}
\label{eq:behPhi}
\end{equation}
(the last equality follows from $\rho^{-1}\Phi\rho=-\Phi$, which
is a consequence of (\ref{eq:defPhi})).
On the other hand, since we took a holomorphic framing the Chern
connection takes the form, in our trivialisation,
$d_H=d+(\partial M_H)M_H^{-1}$
(see for example p. 73 in \cite{GH}), from which we deduce that
$d_{H'}=\rho d_H\rho^{-1}$ and hence $F_{H'}=\rho F_H\rho^{-1}$.
This, together with (\ref{eq:behPhi}), implies that $H'$ is another
solution to (\ref{eq:EHh}). Consequently, $H'':=H+H'$ is also a 
solution to (\ref{eq:EHh}) which satisfies $H''=h_{V}+h_{V^*}$
for some metric $h_V$ (resp. $h_{V^*}$) on $V$ (resp. $V^*$).
To finish the argument, observe that we have an isomorphism
$f:E\to E^*$ defined as $f(u,v)=(v,u)$; now, if $H$ is a solution to
(\ref{eq:EHh}), then both $f^*H$ and
$H^*$ (the latter denotes the natural metric on $E^*$ induced by $H$)
give solutions to the Hermite--Einstein equations for the
pair $(E^*,\Phi^*)$. It follows that $H''+(f^{-1})^*(H'')^*$
is a solution to (\ref{eq:EHh}) which is of the form $h+h^*$
for some metric $h$ on $V$. This metric gives then a solution to
(\ref{eq:HE0}).
\end{pf}

\begin{remark}
Theorem \ref{thm:HK} follows also from the general Hitchin--Kobayashi 
correspondence proved in \cite{BGM}.
\end{remark}

\subsection{Homeomorphism between $\RRR(d)$ and $\MMM(d)$}
\label{sec-gauge-higgs}

The moduli space $\MMM(d)$ of poly\-stable $\Sp(4,\RR)$-Higgs bundles
can be constructed in essentially the same way as that of the moduli
space of ordinary Higgs bundles (see \S 9 of \cite{Si3}).
Alternatively one can do the following: consider the map from 
isomorphism classes of polystable $\Sp(4,\RR)$-Higgs bundles
to isomorphism classes of polystable Higgs bundles
which sends $(V,\varphi)$ to $(E,\Phi)$; this map is finite to one,
so there is a unique structure of reduced scheme on the set of
isomorphism classes of polystable $\Sp(4,\RR)$-Higgs bundles
which is compatible both with the map and with the scheme structure
on the set of isomorphism classes of polystable Higgs bundles
which one gets by looking at the latter as a moduli space.

Following either of these procedures, one ends up in particular
with a structure of topological space on $\MMM(d)$.
We will sketch here a description of this topological structure
on $\MMM(d)$ in gauge-theoretic terms.
This discussion will also be relevant for us in order
to define the  
proper function from which we plan to extract information about the number 
of connected  components of $\MMM(d)$ (a similar study can be seen in  
\S 5 of \cite{Hi1}).

Let $\bV$ be a $C^{\infty}$ complex vector bundle over $X$ of 
rank $2$ and degree $d$. Let us fix a Hermitian metric $h$ on $\bV$,
and denote by $\AAA$ the affine space of Hermitian connections
on $\bV$. Let also $\Omega:=\Omega^{1,0}(M;S^2\bV\oplus S^2\bV^*)$
and let $\GGG$ be the group of Hermitian automorphisms of $\bV$ --- the gauge 
group. Consider the completions of $\AAA$ and
$\Omega$ (resp. $\GGG$) with respect to Sobolev $L^2_1$ (resp. $L^2_2$)
norm, and denote the resulting completions by the same symbols.
Then the quotient $\BBB:=(\AAA\times\Omega)/\GGG$ is a Hausdorff topological
space (the crucial point, Hausdorffness, follows from the existence
of slices of the action of $\GGG$ on $\AAA\times\Omega$).
Finally, we define $\SSS(d)$ to be the set of gauge equivalence
classes $[A,\varphi]\in \BBB$ such that
\begin{equation}
\left\{
\begin{array}{l}
\ov{\partial}_A\varphi=0,\\
F_A+(\beta\beta^*-\gamma^*\gamma)=0, 
\end{array}
\right.
\label{gauge-eq}
\end{equation}
where $\varphi=(\beta,\gamma)$.  
The set $\SSS(d)$ inherits a topology from its
inclusion in $\BBB$. By Theorem \ref{thm:HK},
the points of $\SSS(d)$ are naturally in bijection with
$\MMM(d)$, the moduli space of polystable $\Sp(4,\RR)$-Higgs bundles
of degree $d$ . This bijection maps a gauge equivalence class $[A,\varphi]$
to the pair $(V,\varphi)$, where $V$ is $\bV$ equipped with  the holomorphic 
structure defined by the $(0,1)$ part of the connection $A$.

The following is an $\Sp(4,\RR)$-version of the general correspondence
between complex representations of $\pi$ and Higgs bundles 
(\cite{Hi1,Hi2,Si1,Si2,Do,Co}).

\begin{theorem}
Let $d$ be an integer. There is a homeomorphism $\RRR(d)\cong\MMM(d)$
which, when $d\neq 0$, restricts to a homeomorphism between
the space of irreducible representations in $\RRR(d)$ and the
space of stable Higgs bundles in $\MMM(d)$.
\label{thm:correspondence}
\end{theorem}
\begin{pf}
In fact we prove  that  both $\MMM(d)$ and $\RRR(d)$  are homeomorphic
to $\SSS(d)$.
  To see that $\SSS(d)$ is homeomorphic to
  $\MMM(d)$ we can consider the latter space from the
  complex analytic point of view
  below): consider pairs $(\ov{\partial}_V,\varphi)$,
 where $\ov{\partial}_V$ is  a $\ov{\partial}$-operators on 
the $C^\infty$ vector bundle  $\bV$ underlying $V$
and $\varphi\in\Omega$.  Let $\CCC$ be the set
  of such pairs for which $\varphi$ is holomorphic and the associated
  $\Sp(4,\RR)$-Higgs bundle is polystable.  We can then view
  $\MMM(d)$ as the quotient of $\CCC$ by the
  complex gauge group.  We clearly have an inclusion of the 
space of pairs $(A,\varphi)\in \AAA\times \Omega$ which solve
(\ref{gauge-eq}) into  $\CCC$, which descends to give a continuous map from
  $\SSS(d)$ to $\MMM(d)$. Theorem \ref{thm:HK}
  now shows that this map is in fact a homeomorphism.
 
Suppose that $(E=V\oplus V^\ast,\Phi)$ represents a point in $\MMM(d)$, i.e.
suppose that it is a polystable $\Sp(4,\RR)$-Higgs bundle of degree $d$.
From Theorem \ref{thm:HK}, there is a metric $h$ in $V$ satisfying 
(\ref{eq:HE0}). Rewriting the equations in terms of the Higgs connection
$D=d_A+\theta$, where $A$ is the metric connection and 
$\theta=\Phi+\Phi^\ast$, we see that $D$ is a flat $\Sp(4,\RR)$-connection
and thus defines a point in $\RRR(d)$.
Notice that $d_A$ takes values in $\ulie(2)$, while $\theta$ takes values 
in the orthogonal (w.r.t. the Killing pairing)
complement of $\ulie(2)\subset\splie(4,\RR)$. 
Conversely, by Corlette's theorem \cite{Co}, every representation
in $\RRR(d)$ arises in this way. The fact that this correspondence gives
a homeomorphism follows by the same argument as the one given in 
\cite{Si3} for ordinary Higgs bundles.
When $d\neq 0$, the solution to (\ref{gauge-eq}) is irreducible if and only 
if the corresponding $\Sp(4,\RR)$-Higgs bundle is stable, hence  the 
corresponding element 
in $\RRR(d)$ is irreducible. When $d=0$ the solution to (\ref{gauge-eq}) in 
 a polystable $\Sp(4,\RR)$-Higgs
bundle may be actually $\Sp(4,\RR)$-irreducible. The reason for this lies in 
the fact 
that we have used the standard  stability  of $(E,\Phi)$ as the stability
criterium for the $\Sp(4,\RR)$-Higgs bundle. There is, however,  a notion
of stability of the $\Sp(4,\RR)$-Higgs bundle in its own rigth, i.e. 
without
using that of the corresponding Higgs bundle $(E,\Phi)$
(\cite{BGM}).
It turns out that these two notions are equivalent when $d\neq 0$, 
but when $d=0$, stability of the Higgs bundle $(V,\varphi)$ is only 
equivalent to polystability.
This is, however, no so important for us since in Theorem \ref{higgs-moduli}, 
where this is used, we assume that $d\neq 0$.
\end{pf}

\subsection{A proper function on $\MMM(d)$}

We follow the ideas of \cite{Go2,Hi2}, which reduce the proof of
the connectedness of $\MMM(d)$ to proving connectedness of a smaller
subspace $\NNN(d)\subset\MMM(d)$.
Let us briefly recall how this goes.

We  define, for any $[A,\varphi]\in\SSS(d)\cong\MMM(d)$,
$$f([A,\varphi]):=\|\beta\|_{L^2}^2+\|\gamma\|_{L^2}^2,$$
where $\varphi=(\beta,\gamma)$.
This expression is gauge invariant and hence descends
to give a map
$$f:\MMM(d)\to\RR.$$
One can prove that $f$ is proper, essentially by using Uhlenbeck's
compactness theorem (see \cite{Hi1}). So any connected
component of $\MMM(d)$ contains a local minimum of $f$ and
hence, we have the following.
\begin{prop}
Let $\NNN(d)\subset\MMM(d)$ be the local minima of $f$. Then
$\MMM(d)$ is connected if $\NNN(d)$ is connected.
\end{prop}
By duality we have  $\MMM(d)\cong
\MMM(-d)$, so it suffices to consider the case $d<0$.
The key result characterizing $\NNN(d)$ is the following.
\begin{prop}[Gothen, \cite{Go2}]
\label{minima}
Suppose that $d$ satisfies $-2(g-1)<d<0$. Then
$\NNN(d)$ consists of the
classes $[A,\varphi]$ such that $\gamma=0$.
\end{prop}

The following sections are devoted to proving that  
if $-2(g-1)<d<0$, the subspace
$\NNN(d)$ is connected and contains stable $\Sp(4,\RR)$-Higgs bundles, 
thus finishing the proof of Theorem \ref{higgs-moduli}, and hence
the main theorems stated in the introduction.
To study $\NNN(d)$, we take a more general point of view and consider
a moduli problem that depends on a real parameter $\alpha$.
 We can then  identify
 $\NNN(d)$ with the  moduli space for $\alpha=0$.

\section{Quadratic pairs}\label{sec-pairs}

 \subsection{} 

A {\bf quadratic pair} $(V,\beta)$ on $X$ consists by definition of 
a holomorphic vector bundle $V$ on $X$ of rank $2$ and a holomorphic section
$\beta\in H^0(K\otimes S^2V)$, where $K$ is the canonical
bundle of $X$. The degree of $(V,\beta)$ is by definition the degree of $V$.
Let $\alpha\in\RR$. We say that $(V,\beta)$ is
{\bf $\alpha$-polystable} if, denoting by $d$ the degree of $(V,\beta)$,
 \begin{enumerate}
 \item[(1)] $d/2\leq\alpha$, and if $\beta=0$ then $\alpha=d/2$,
 \item[(2)] for any subbundle $L\subset V$
 \begin{enumerate}
 \item if $\beta\in H^0(K\otimes S^2L)$ then $\deg L \leq d-\alpha$,
and in case there is equality there is a splitting $V=L\oplus L'$;
 \label{condicio1}
 \item if $\beta\in H^0(K\otimes L\otimes V)$ then
$\deg L \leq d/2$, and if there is equality then we have
a splitting $V=L\oplus L'$ in such a way that
$\beta\in H^0(K\otimes L\otimes L')$; 
 \label{condicio2}
 \item $\deg L \leq \alpha$ in any case, and if there is equality
then there is a splitting $V=L\oplus L'$.
 \label{condicio3}
 \end{enumerate}
 \end{enumerate}
We say that $(V,\beta)$ is {\bf $\alpha$-stable} if equality
never occurs in the inequalities required for polystability.

For any rational value of  $\alpha$ there exists an algebraic coarse
moduli space $\NNN_{\alpha}(d)$
for the moduli problem of families of $\alpha$-polystable quadratic
pairs over $X$ of degree $d$. As in the case of the moduli space
of symplectic Higgs bundles, we will not describe the algebraic construction
of $\NNN_{\alpha}(d)$, since for our purposes it suffices to have
a topological description of it.
To construct the topological space underlying $\NNN_{\alpha}(d)$ one uses
the same strategy as in Section  \ref{sec-gauge-higgs}. Namely, to use 
a version of the Hitchin--Kobayashi correspondence for
$\alpha$-polystable quadratic pairs (analogous to Theorem
\ref{thm:HK}), a particular case of the more general correspondence proved
in  \cite{BGM},  and then use gauge theory.

It is clear that a quadratic pair is a special case of 
$\Sp(4,\RR)$-Higgs bundle. The next proposition  says
that when $\alpha=0$ the two possible notions of (poly)stability coincide.

\begin{prop}
\label{prop:stabd}
Let $-2(g-2)<d<0$ be an integer, and let $(V,\beta)$ be a quadratic
pair of degree $d$. Then, $(V,\beta)$ is a polystable (resp. stable)
$\Sp(4,\RR)$-Higgs bundle if and only if it is a $0$-polystable
(resp. $0$-stable) quadratic pair. In other words, $\NNN(d)=\NNN_0(d)$.
\end{prop}
\begin{pf}
That polystability as $\Sp(4,\RR)$-Higgs bundle implies $0$-polystability
as quadratic pair is obvious. To prove the converse, 
let $p:V\oplus V^*\to V$ and $q:V\oplus V^*\to V^*$ denote the
projections. Assume that $(V,\beta)$ is $0$-polystable as quadratic
pair, and consider some subbundle $W\subset V\oplus V^*$ satisfying
$\Phi(W)\subset K\otimes W$. We need to check that $\deg W\leq 0$, and
that equality implies splitting. We will prove the first claim, since
the second one follows from very similar argument. 
Define $A=p(W)$, $B=q(W)$, $A'=W\cap V$ and $B'=W\cap V^*$.
Using $0$-polystability it is straightforward to check that
\begin{equation}
\deg(A^{\perp}+B)+\deg(A+B^{\perp})\leq 0
\label{eq:54}
\end{equation}
(simply consider the subbundle $A+B^{\perp}\subset V$).
One computes
\begin{align*}
\deg (A+B^{\perp}) &= \deg A+\deg B-\deg(A^{\perp}+B) \\
\deg (A^{\perp}+B) &= \deg A+\deg B-\deg(A+B^{\perp}).
\end{align*}
Adding up and taking (\ref{eq:54}) into account we get 
$\deg A+\deg B\leq 0$. A similar argument implies that
$\deg A'+\deg B'\leq 0$. Finally, combining these two
inequalities with the exact sequences
$$0\to B'\to W\to A\to 0
\qquad\text{ and }\qquad
0\to A'\to W\to B\to 0$$
we get the required inequality $\deg W\leq 0$.
\end{pf}

\subsection{Familes of quadratic pairs and moduli space}
\label{ss:families}

If $U$ is a scheme of finite type (resp. a complex manifold) we define
a {\bf family of quadratic pairs} on $X$ parametrized by $U$ to be a pair
$(\VVV,\Beta)$, where $\VVV$ is a rank $2$ vector bundle on $U\times X$
and $\Beta\in H^0(\pi_X^*K\otimes S^2\VVV)$ is an algebraic
(resp. holomorphic) section
(here $\pi_X:U\times X\to X$ denotes the projection).

The following lemma states
that stability for algebraic quadratic pairs is an open condition in
the Zariski topology, as one would naturaly expect.
 
\begin{lemma}
Fix a real number $\alpha$.
Let $U$ be a scheme of finite type and let $(\VVV,\Beta)$
be an algebraic family of quadratic pairs on $X$ parametrized by $U$.
For any $u\in U$ let $(\VVV_u,\Beta_u)$ denote the
restriction of $(\VVV,\Beta)$ to
$X\times\{u\}$. The set
$$U_s(\VVV,\Beta)=
\{u\in U\mid (\VVV_s,\Beta_s)\text{ is $\alpha$-stable }\}$$
is a Zariski open subset of $U$.
\label{openness}
\end{lemma}
\begin{pf}
Fix  $d\in\ZZ$ and let $J:=\Jac_dX$.
Let $\PPP\to J\times X$ be the Poincar\'e bundle.
Define $\FFF:=(\pi_U)_*
(\pi_{U\times X}^*\VVV\otimes \pi_{J\times X}^*\PPP^\ast),$
where $\pi_U:U\times J\times X\to U$ is the projection
and the other maps $\pi_{?}$ are the obvious analogs.
Let $\LLL\to J\times X$ be a line bundle of high enough degree
such that the natural map
$\pi_U^*{\pi_U}_*(\pi_{J\times X}^*\PPP^{\ast}\otimes\LLL)\to
\pi_{J\times X}^*\PPP^{\ast}\otimes\LLL$
is surjective,
${\mathbf R}^j{\pi_U}_*(\pi_{J\times X}^*\PPP^{\ast}\otimes \LLL)=0$
for any $j>0$, and the same thing applies for
$\pi_{U\times X}^*\VVV\otimes \LLL$, $\pi_{U\times X}^*\VVV^{\ast}\otimes \LLL$
and $\pi_X^*K\otimes\LLL^2$.
Denote by $\III$ the image sheaf of the canonical map
$\FFF\otimes{\pi_U}_*(\pi_{J\times X}^*\PPP^{\ast}\otimes\LLL)\to
{\pi_U}_*(\pi_{U\times X}^*\VVV\otimes\LLL)$, and let
$\GGG:=\III^{\perp}\subset 
{\pi_U}_*(\pi_{U\times X}^*\VVV^{\ast}\otimes\LLL)$
be the orthogonal subsheaf.
Denote by $Q:\VVV^{\ast}\to K$ be the quadratic map induced by $\Beta$
(of course this is not a morphism of vector bundles!).
The map $Q$ induces another map
$\QQQ:{\pi_U}_*(\pi_{U\times X}^*\VVV^{\ast}\otimes\LLL)
\to{\pi_U}_*(\pi_{X}^*K\otimes\LLL^2)$. Finally, let 
$\HHH:=\GGG\cap\QQQ^{-1}(0)$. This is a coherent sheaf on $U$.
Let $U_d$ be the support of $\HHH$ with the reduced scheme structure.
This is by definition a closed subscheme of $U$ with the property that
a closed point $u\in U$ belongs to $U_d$ if and only if there
is a line subbundle $L\subset \VVV_u$ with degree $\deg L\geq d$ and
such that $B_u\in H^0(L\otimes\VVV_u\otimes K)$. Similarly, one proves
that the subscheme $U_d'\subset U$ whose closed points are the $u\in U$
such that there
is a line subbundle $L\subset \VVV_u$ with degree $\deg L\geq d$ 
satisfying $B_u\in H^0(L^2\otimes K)$ is closed in $U$, and that
the subscheme $U_d''\subset U$ whose closed points are the $u\in U$
such that there
is a line subbundle $L\subset \VVV_u$ with degree $\deg L\geq d$
is also closed in $U$.
\end{pf}

\begin{corollary}
Let $U$ be an irreducible
algebraic manifold and let $(\VVV,\Beta)$ be a (holomorphic)
family of quadratic pairs on $X$ parametrized by $U$. The set
$$U_s(\VVV,\Beta)=\{u\in U\mid (\VVV_s,\Beta_s)\text{ is $\alpha$-stable }\}$$
is connected.
\label{cor:stableconn}
\end{corollary}
\begin{pf}
By Serre's GAGA the holomorphic family $(\VVV,\Beta)$ is induced
by an algebraic family $(\VVV^{\alg},\Beta^{\alg})$. We then have
$U_s(\VVV,\Beta)=U_s(\VVV^{\alg},\Beta^{\alg})$, and the latter is
connected in the standard topology, since it is a Zariski open
subset of the irreducible manifold $U$.
\end{pf}

The following continuous version of the algebraic notion of family
will also be relevant for us. We say that a quadratic pair $(V,\beta)$
{\bf can be (continuously) deformed to} another
pair $(V',\beta')$ if and only if $V$ and $V'$ are isomorphic as smooth
vector bundles and, denoting by $\bV$ the $C^{\infty}$ complex vector bundle
on $X$ underlying both $V$ and $V'$ and by $\ov{\partial}_V$ (resp.
$\ov{\partial}_{V'}$) the $\ov{\partial}$-operator giving rise to the
holomorphic structure of $V$ (resp. $V'$), there are continuous maps
$D:[0,1]\to\Omega^{0,1}(\End \bV)$ and $B:[0,1]\to\Omega^0(S^2\bV\otimes K)$ 
such that:
\begin{itemize}
\item[(i)] $D(0)=0$ and $D(1)=\ov{\partial}_{V'}-\ov{\partial}_V$;
\item[(ii)] $B(0)=\beta$ and $B(1)=\beta'$;
\item[(iii)] for any $t\in[0,1]$ we have
$(\ov{\partial}_V+D(t))B(t)=0$.
\end{itemize}
Denote by $(V_t,\beta_t)$ the quadratic pair defined by 
$(\ov{\partial}_V+D(t),B(t))$. We say that the continuous deformation
$(D,B)$ {\bf goes through $\alpha$-(poly)stable pairs} if for any $t$ the pair
$(V_t,\beta_t)$ is $\alpha$-(poly)stable.

 \subsection{Geometry of quadratic pairs}

 \subsubsection{Quadratic forms in $\CC^2$}
 \label{sect:quad}
 Let $\VV=\CC^2$, and let $x,y\in \VV$ be a basis. Recall
 that the {\bf discriminant} of $f=ax^2+bxy+cy^2\in S^2\VV$ is
 by definition 
 $\Delta_f:=(b^2-4ac)(x\wedge y)^2\in (\Lambda^2\VV)^2.$
 It is clear that $\Delta:S^2\VV\to (\Lambda^2\VV)^2$
 is a $\GL(\VV)$ equivariant map.
 Define the following subsets of $S^2\VV$:
 $\OOO_0 = \{0\}$,
 $\OOO_1 = \{f\in S^2\VV\mid f\neq 0,\ \Delta_f=0\}$ 
 and $\OOO_2 = \{f\in S^2\VV\mid f\neq 0,\ \Delta_{f}\neq 0\}$.

 \begin{lemma}
 {\rm (i)} An element $0\neq f\in S^2\VV$ belongs to
 $\OOO_1$ if and only if there exists $0\neq x\in \VV$ so that
 $f=x^2$, and in this case the span $\CC x\subset \VV$ only depends
 on $f$. 
 {\rm (ii)} $f$ belongs to $\OOO_2$ if and only if there
 exist linearly independent elements $x,y\in \VV$ so that
 $f=xy$, and in this case the set of lines $\{\CC x,\CC y\}$ only
 depends on $f$.
 {\rm (iii)} The orbits of the action of $\GL(\VV)$ on $S^2\VV$ are 
 $\OOO_0,\OOO_1,\OOO_2$.
 \label{lemma:orbit}
 \end{lemma}
 \begin{pf}
 (i) and (ii) are easy, and (iii) follows from them.
 \end{pf}

 One readily sees (for example, by computing the dimension of
 the stabilizers of $\GL(\VV)$ acting on $S^2\VV$)
 that $\OOO_2\subset S^2\VV$ is open. This implies the following.

 \begin{lemma}
 Let $f\in S^2\VV$, and let $\rho_f:\End(\VV)\to T_f(S^2\VV)\cong S^2\VV$
 be the map given by the infinitesimal action of $\End(\VV)=\Lie\GL(\VV)$
 on $S^2\VV$. If $f\in\OOO_2$, then $\rho_f$ is onto.
 \label{lemma:orbites}
 \end{lemma}

 The following lemma is straightforward.

 \begin{lemma}
 Let $f\in S^2\VV$, and let $\rho_f:\End(\VV)\to T_f(S^2\VV)\cong S^2\VV$
 be as above. 
 If $f=xy\in\OOO_2$ then, in the basis $\la x,y\ra$ of $\VV$, 
 $$\Ker \rho_f=\left\{\left(\begin{array}{cc}
 \lambda & 0 \\ 0 & -\lambda\end{array}\right)\mid \lambda\in\CC\right\}.$$
 %
 \label{lemma:kerrho}
 \end{lemma}

 \begin{lemma}
 If $f,g\in S^2\VV$ are linearly independent, then for a generic
 value of $\lambda\in\CC$ $f+\lambda g$ belongs to $\OOO_2$.
 \label{lemma:nolinia}
 \end{lemma}
 \begin{pf}
 We have to check that $\OOO_1$ does not contain any plane.
 This is equivalent to proving that the subset 
 $\PP\OOO_1\subset\PP(S^2\VV)$
 induced by $\OOO_1$ does not cointain any line. But
 $\PP\OOO_1$ is given by the equation 
 $4ac=b^2$, hence is a nondegenerate conic and does not
 contain any line.
 \end{pf}

 \subsubsection{Discriminant of quadratic pairs}
 \label{casosDelta}
 Extending fibrewise the definition of the discriminant
we get a quadratic map 
 $\Delta:K\otimes S^2V\to K^2\otimes (\Lambda ^2V)^2.$
 For any section $\beta:\cO\to K\otimes S^2V$
 we denote $\Delta\circ\beta$ by $\Delta_{\beta}$.
 Note that we always have $\beta^{-1}(0)\subset\Delta_{\beta}^{-1}(0)$.
 It follows from Lemma \ref{lemma:orbit} that, given a pair $(V,\beta)$,
 we can distinguish four possibilities.
 \begin{enumerate}
 \item[(1)] $\beta=0$.
 \item[(2)] $\beta\neq 0$ and $\Delta_{\beta}=0$. Then there exists a line
 subbundle $L\subset V$ so that $\beta\in H^0(K\otimes L^2)$.
 \item[(3)] $\Delta_{\beta}\neq 0$ and there exists a square root
 $\Delta^{1/2}(\beta)\in H^0(K\otimes\Lambda^2V)$; then there exist
 two different line subbundles $L_1,L_2\subset V$ so that 
 $\beta\in H^0(K\otimes L_1\otimes L_2)$.
 \item[(4)] There is no square root of $\Delta_{\beta}$; then there exists
 no line subbundle $L\subset V$ so that $\beta\in H^0(K\otimes L\otimes
 V)$.
 \end{enumerate}
 To understand the difference between (3) and (4), observe that by
 Lemma \ref{lemma:orbit} for any $x\in X$ such that $\Delta_{\beta}(x)\neq 0$
 one has two different lines $L_{x,1},L_{x,2}\subset V_x$. As $x$
 moves around $X\setminus\Delta_{\beta}^{-1}(0)$ these two lines 
 give rise to a pair of line subbundles of $V$, unless the monodromy
 around points
 in $\Delta_{\beta}^{-1}(0)$ interchanges the lines. Now, the existence
 of a square root of $\Delta_{\beta}$ is equivalent to the
 triviality of the action of the monodromy on the pair of lines.

 \begin{lemma}
 \label{casD12}
 Suppose that a pair $(V,\beta)$ defined on $X$ 
 satisfies $\Delta=\Delta_{\beta}\neq 0$ and
 $\Delta^{1/2}$ exists. Let $L_1,L_2\subset V$ be the 
 two different line bundles so that
 $\beta\in H^0(K\otimes L_1\otimes L_2)$.
 Then:
 {\rm (i)} the pair $\{L_1,L_2\}$ is uniquely determined by $V$ and $\beta$;
 {\rm (ii)} there is an exact sequence of sheaves 
 $$0\to L_1\oplus L_2\to V\to \cO_T\to 0;$$
 here $T$ is the divisor in $X$ defined as
 $T=\sum_{x\in X}(r(x)/2-z(x))x$,
 where $r(x)$ (resp. $z(x)$) is the vanishing order of $\Delta$
 (resp. $\beta$) at $x$.

 \end{lemma}
 \begin{pf}
 (i) is true because the induced quadratic map $\gamma:V^*\to K$ vanishes
 exactly at $L_1^{\perp}\cup L_2^{\perp}$; (ii) follows from the next
 lemma.
 \end{pf}

 \begin{lemma}
 \label{estudilocal}
 Let $(V,\beta)$ be a pair defined on some disk $\DD\subset X$,
 and assume that $\Delta=\Delta_{\beta}\neq 0$, $\Delta^{-1}(0)=\{0\}$, and
 there is a square root $\Delta^{1/2}:\DD\to\CC$ of $\Delta$. 
 Let $r$ (resp. $z$) be the vanishing order of $\Delta$ (resp. $\beta$)
 at $0$. Let $\theta$ be a coordinate in $\DD$.
 Then:
 {\rm (i)} we have $r/2\geq z$;
 {\rm (ii)} let $L_1,L_2\subset V$ be the line subbundles so that
 $\beta\in H^0(K\otimes L_1\otimes L_2)$; 
 we have an exact sequence of sheaves:
 $0\to L_1\oplus L_2\to V\to \cO/\theta^{r/2-z}\cO\to 0$.
 \end{lemma}
 \begin{pf}
 Pick a trivialisation $V\cong\DD\times\CC\la x,y\ra$ so that
 $L_1=\DD\times\CC\la x\ra$. Then we can write $\beta=x(ax+by)$,
 where $a,b$ are holomorphic functions on $\DD$, and we have
 $\Delta=b^2$. We then have: $a=\theta^z(a_0+\theta a_1)$ and
 $b=\theta^{r/2}(b_0+\theta b_1)$, where $a_0,b_0\in\CC^{\times}$ and
 $a_1,b_1$ are holomorphic. To prove (i), observe that $\theta^z$
 divides $b$. (ii) follows from $L_2=\CC\la ax+by\ra$.
 \end{pf}

 \subsection{Local study of moduli space}

 Take an isomorphism class $[V,\beta]\in\NNN_{\alpha}(d)$.
 By a result of Biswas and Ramanan \cite{BiRa} the Zariski tangent space 
 of $\NNN_{\alpha}(d)$ at
 $[V,\beta]$ is given by the first hypercohomology group 
 of the following $2$-term complex:
 $$\CCC_{V,\beta}:\End V\ni \psi \mapsto \rho_{V,\beta}(\psi)\in K\otimes 
S^2 V,$$
 where $\rho_{V,\beta}$ is the map induced fibrewise by the infinitesimal
 action of $\End(\CC^2)$ on $S^2\CC^2$ along $\beta\subset K\otimes
 S^2V$.
 Applying hypercohomology to this exact sequence of complexes
 $$0\to K\otimes S^2V[-1]\to\CCC_{V,\beta}\to\End V[0]\to 0$$
 (as usual, if $E$ is a sheaf, $E[d]$ denotes the complex whose
 only nonzero term is a copy of $E$ in the position $-d$)
 we get the following long exact sequence:
 \begin{align}
 0 &\to \HH^0(\CCC_{V,\beta})\to H^0(\End V)\to H^0(K\otimes S^2V)\to 
 \label{eq:longHH}\\
 &\to\HH^1(\CCC_{V,\beta})\to H^1(\End V)\to H^1(K\otimes S^2V)\to
 \HH^2(\CCC_{V,\beta})\to 0. \notag
 \end{align}

 \begin{lemma}
 Suppose that $\HH^0(\CCC_{V,\beta})=\HH^2(\CCC_{V,\beta})=0$.
 Then $\dim\HH^1(\CCC_{V,\beta})=7(g-1)+3d.$
 \end{lemma}
 \begin{pf}
 Using (\ref{eq:longHH}) and Riemann--Roch
 we have 
 \begin{align*}
 \dim\HH^1(\CCC_{V,\beta})&=-\chi(\CCC_{V,\beta})=
 -\chi(\End V)+\chi(K\otimes S^2V) \\
 &=-4(1-g)+(3(2g-2)+3d+3(1-g))=7(g-1)+3d.
 \end{align*}
 \end{pf}

 Biswas and Ramanan also prove that the points $[V,\beta]\in\NNN_{\alpha}(d)$
 at which $\HH^0=\HH^2=0$ are smooth. The preceeding lemma
 (together with the identification $T_{V,\beta}\NNN_{\alpha}(d)\cong \HH^1$
 given by Biswas and Ramanan) implies that the dimension of $\NNN_{\alpha}(d)$
 at these points is $7(g-1)+3d$.

 \subsection{Vanishing of $\HH^0$}
  
 \begin{lemma}
 Suppose that $(V,\beta)$ is $\alpha$-polystable for some
 $\alpha$ and $\Delta_{\beta}\neq 0$.
 Then we have $\HH^0(\CCC_{V,\beta})=0$ unless there is
 a splitting $V=L_1\oplus L_2$ in line bundles with
 $\deg L_1=\deg L_2=d/2$ and $\beta\in H^0(K\otimes L_1\otimes L_2)$.
 \label{lemma:vanishH0}
 \end{lemma}
 \begin{pf}
 Thanks to (\ref{eq:longHH}) the vanishing of $\HH^0$ is equivalent
 to the injectivity of the map $H^0(\End V)\to H^0(K\otimes S^2V)$.
 Suppose the map is not injective.
 Then there is some nonzero $s\in H^0(\End V)$ so that
 $\rho_{V,\beta}(s)=0$.

 Since $\Delta_{\beta}\neq 0$, Lemma \ref{lemma:kerrho} implies that 
 if $x\in X$ and $s(x)\neq 0$ then $s(x)$ has two different
 eigenvalues. On the other hand, since the characteristic polynomial
 of $s$ has holomorphic coefficients in $X$, hence constant
 coefficients, the eigenvalues of $s(x)$ are constant as $x$ varies
 along $X$. Hence we can split $V=L_1\oplus L_2$ in eigensubbundles 
 of $s$. Using again Lemma \ref{lemma:kerrho}, we know that 
 $\beta\in H^0(K\otimes L_1\otimes L_2)$. Applying the semistability condition
 we get $\deg L_1\leq d/2$, $\deg L_2\leq d/2$, and hence
 $\deg L_1=\deg L_2=d/2$.
  \end{pf}

 \subsection{Vanishing of $\HH^2$}

 \begin{lemma}
 Assume that $\Delta_{\beta}\neq 0$. Then $\HH^2(\CCC_{V,\beta})=0$.
 \end{lemma}
 \begin{pf}
 It suffices to show that the map $H^1(\End V)\to H^1(K\otimes S^2V)$
 in (\ref{eq:longHH}) is onto (this map is the one induced by
 $\rho=\rho_{V,\beta}$ in cohomology). By Serre duality, this is equivalent
 to the injectivity of the map
 $H^0(S^2V^*)\to H^0(K\otimes \End V^*)$
 induced by $\tau_K(\rho^*):S^2V^*\to K\otimes\End V^*$.
 But it follows from Lemma \ref{lemma:orbites} that the map
 of sheaves $\rho:\End V\to K\otimes S^2V$ is onto, so 
 the map $\rho^*:K^*\otimes S^2V^*\to\End V^*$ is injective.
 Then, by Lemma \ref{lemma:ABL} below,
 $\tau_K(\rho^*):S^2V^*\to K\otimes\End V^*$
 is also injective.
 \end{pf}

 \begin{lemma}
 Let $L$ be a line bundle on $X$ and let $A,B$ be coherent
 sheaves on $X$. Let 
 $\tau_L:\Hom(L^*\otimes A,B)\to\Hom(A,L\otimes B)$
 be the standard isomorphism and let $f\in\Hom(L^*\otimes A,B)$.
 Then $\Ker \tau_L(f)=L \otimes \Ker f.$
 \label{lemma:ABL}
 \end{lemma}
 \begin{pf}
 The sheaf  $L$ is flat since it is locally free, and hence the functor
 $L\otimes\cdot:\Coh(X)\to\Coh(X)$ is exact. But $\tau_L$ is the
 map induced by $L\otimes\cdot$ on morphisms in $\Coh(X)$.
 \end{pf}

 \section{Proof of Theorem \ref{thm:conn}} 
\label{sec-conn}

 We briefly sketch our
strategy. We first classify the pairs $(V,\beta)$ in three
types as follows:
 \begin{itemize}
 \item {\bf type A}: pairs $(V,\beta)$ with $\Delta_{\beta}=0$;
 \item {\bf type B}: pairs $(V,\beta)$ with $\Delta_{\beta}\neq 0$ and
 $\HH^0(\CCC_{V,\beta})\neq 0$;
 \item {\bf type C}: pairs $(V,\beta)$ with $\Delta_{\beta}\neq 0$ and
 $\HH^0(\CCC_{V,\beta})= 0$.
 \end{itemize}
In section \ref{ss:typeA} we prove that any $\alpha$-polystable pair 
 of type A can be deformed to a $\alpha$-polystable
 pair of type B or C, and
in section \ref{ss:typeB} we prove that any $\alpha$-polystable pair
of type B can be deformed to a $\alpha$-polystable pair of type C (in both
cases, and everywhere below, we mean {\it continuous} deformations).
Finally, in section \ref{ss:typeC} we prove that the subset of
$\NNN_{\alpha}(d)$ consisting of pairs of type $C$ is indeed connected.
In Subsection \ref{ss:existencia} we prove that $\NNN_{\alpha}(d)$ contains
at least one stable object, thus proving the last claim in Theorem
\ref{thm:conn}.

 \subsection{Pairs of type A}
 \label{ss:typeA}

 In this section we study the $\alpha$-polystable pairs $(V,\beta)$ for which
 $\Delta_{\beta}=0$.
 By Subsection \ref{casosDelta}, for any such pair there is a line subbundle
 $L\subset V$ so that $\beta\in H^0(K\otimes L^2)$.
 Let $l$ be the degree of $L$. Since $0\neq\beta\in H^0(K\otimes L^2)$, we have
 $\deg K\otimes L^2=2g-2+2l\geq 0$. Also,
 if $(V,\beta)$ is $\alpha$-polystable for some $\alpha\geq 0$ then we must have
 $l\leq d$. So in this section we assume that $d\geq -(g-1)$.

 \begin{lemma}
 There exists a manifold $S_{d,l}$ and a family of pairs
 $(\VVV,\Beta)$ on $X$ parametrized by $S_{d,l}$ such that 
 any pair $(V,\beta)$ with 
 $\deg V=d$, $\beta\in H^0(K\otimes L^2)$ and $\deg L=l$,
 is isomorphic to $(\VVV|_{\{s\}\times X},\Beta|_{\{s\}\times X})$
 for at least one $s\in S_{d,l}$.
 \end{lemma}
 \begin{pf}
 Let $j=2g-2+2l$, and let $\mu:S^jX\to\Jac_j(X)$ be the map which
 sends an effective divisor $D$ to the bundle $\cO(D)$. 
 Let also $q:\Pic^l(X)\to\Jac_j(X)$ the map defined as
 $q(L)=K\otimes L^2$. Finally, let 
 $\Sigma=S^jX\times_{\Jac_j(X)}\Jac_l(X).$
 Since $\mu$ is \'etale it follows that $\Sigma$ is a manifold.
 Let $\lambda:\Sigma\to\Jac_l(X)$ be the projection.
 We define $S_{d,l}$ to be the vector bundle over 
 $\Sigma \times \Jac_{d-l}(X)$ whose fibre over $(\sigma,L')$
 is $\Ext^1(L',\lambda(\sigma))$ (this vector space has constant
dimension because for any $(L,L')\in\Jac_l(X)\times\Jac_{d-l}(X)$
 we have $H^0(X;L\otimes {L'}^*)=0$, since
 $\deg L\otimes {L'}^*=2l-d<0$).
 It is clear how to construct the family $(\VVV,\Beta)$.
 \end{pf}

 \begin{lemma}
 Suppose that $-(g-1)<d\leq 0$. There is a nonempty
 dense Zariski open subset
 $S_{d,l}^*\subset S_{d,l}$ such that for any $s\in S_{d,l}^*$
 the pair $(\VVV|_{\{s\}\times X},\Beta|_{\{s\}\times X})$
 can be deformed to a pair with $\Delta\neq 0$
 through a path of $\alpha$-stable pairs. 
 \label{lemma:casfacil0}
 \end{lemma}
 \begin{pf}
 Let us prove that generically 
 $h^0(K\otimes L^2)<h^0(K\otimes S^2V)$.
 This means that there is some section $\beta'$ of $K\otimes S^2$
 which is not entirely contained in $K\otimes L^2$.
 Then Lemma \ref{lemma:nolinia} implies that for 
 generic $\lambda\in\CC$ the section $\beta(\lambda)=\beta+\lambda\beta'$
 satisfies $\Delta_{\beta(\lambda)}\neq 0$, thus proving the claim.

 First of all, we have
 $\chi(K\otimes L^2)<\chi(K\otimes S^2V).$
 Indeed, $\chi(K\otimes S^2V)=3(g-1+d)$ and
 $\chi(K\otimes L^2)=g-1+2l$. Now,
 $$g-1+2l<3(g-1+d) \Longleftrightarrow 0<(2(g-1)+d)+(2d-2l),$$
 and the right hand side follows from $2(g-1)+d>0$ (because
 $d\geq -(g-1)$) and $2d-2l\geq 0$ (by stability).
 On the other hand, since $-(g-1)<d$ we have
 \begin{equation}
 \label{eq:casigual1}
 h^0(K\otimes S^2V)\geq\chi(K\otimes S^2V)\geq 3.
 \end{equation}
 Now it follows from standard Brill--Noether theory 
 \cite{ACGH} that
 for a generic $L$ in $\lambda(\Sigma)\subset\Jac_l(X)$ we have 
 $h^0(K\otimes L^2)=\max\{1,\chi(K\otimes L^2)\}.$
 If $h^0(K\otimes L^2)=\chi(K\otimes L^2)$ then we are done,
 since we have
 $\chi(K\otimes L^2)<\chi(K\otimes S^2V)\leq h^0(K\otimes S^2).$
 If instead $h^0(K\otimes L^2)=1$, then we simply use
 (\ref{eq:casigual1}).
 \end{pf}

\begin{lemma}
\label{lemma:casfacil}
Suppose that $-(g-1)<d<0$. Any $\alpha$-stable element in the family
$(\VVV,\Beta)$ can be deformed to a pair with $\Delta\neq 0$ through
a path of $\alpha$-stable pairs.
\end{lemma}
\begin{pf}
Let $S_{d,l}^s\subset S_{d,l}$ be the set of points
$s\in S_{d,l}^s$ for which $(\VVV_{\{s\}\times X},
\Beta_{\{s\}\times X})$ is $\alpha$-stable. By Corollary \ref{cor:stableconn},
$S_{d,l}^s$ is connected. If $S_{d,l}^s$ is empty, then there
is nothing to prove.
Otherwise, by Lemma \ref{lemma:casfacil0} the intersection
$S_{d,l}^s\cap S_{d,l}^*$ is nonempty, and the result follows.
\end{pf}

 \begin{lemma}
 A $\alpha$-polystable element in the family
 $(\VVV,\Beta)$ parametrized by $S_{-(g-1),-(g-1)}$ can be
 deformed to a pair with $\Delta\neq 0$
 through a path of $\alpha$-stable pairs.  
 \end{lemma}
 \begin{pf}
By definition, a $\alpha$-polystable element
$(V,\beta)$ appearing in the family $S_{-(g-1),-(g-1)}$ must split
$V=L\oplus L'$, where $L^2=K$ and $\deg L'=0$. 
Now, $S^2=L^2\oplus L\otimes L'\oplus {L'}^2$ and
$h^0(K\otimes {L'}^2)\geq\chi(K\otimes {L'}^2)=g-1$, so
$h^0(K\otimes L^2)<h^0(K\otimes S^2V)$. The result now
follows from the same argument as in the preceeding lemma.
 \end{pf}

 \subsection{Pairs of type B}
 \label{ss:typeB}

 By Lemma \ref{lemma:vanishH0} the pairs $(V,\beta)$ 
 which are $\alpha$-polystable for some value of $\alpha$ and
 so that $\Delta_{\beta}\neq 0$ and $\HH^0(\CCC_{V,\beta})\neq 0$
 are precisely those of the form
 \begin{equation}
 V=L_1\oplus L_2, \qquad
 \deg L_1=\deg L_2=d/2, \qquad
 \beta\in H^0(K\otimes L_1\otimes L_2).
 \label{eq:aquestsparells}
 \end{equation}
 In this section we assume that $\delta=d/2$ is an integer.

 \begin{lemma}
 Any pair $(V,\beta)$ of the form (\ref{eq:aquestsparells}) is 
 $\alpha$-polystable for any $\alpha\geq 0$.
 \label{totestable}
 \end{lemma}
 \begin{pf}
 If $L\subset L_1\oplus L_2$ is a line subbundle, then
 at least for one value of $i$ the projection $L\to L_i$ is
 nonzero, from which it follows that $\deg L\leq \deg L_i=d/2<0$.
 We cannot have $\beta\in H^0(K\otimes L^2)$,
 so $\deg L\leq d/2<0$ implies $\alpha$-polystability for any $\alpha\geq 0$.
 \end{pf}

 \begin{lemma}
 There exists an irreducible manifold $S_d$ an a family of pairs
 $(\VVV,\Beta)$ on $X$ parametrized by $S_d$ such that
 any pair $(V,\beta)$ 
 of the form (\ref{eq:aquestsparells})
 is isomorphic to
 $(\VVV|_{\{s\}\times X},\Beta|_{\{s\}\times X})$
 for at least one $s\in S_d$.
 \label{lemma:irrfamdelta}
 \end{lemma}
 \begin{pf}
 Define $S_d=S^{2g-2+d}X\times\Jac_\delta(X)$.
 Let $\pi:S_d\times X
 \to\Jac_{\delta}(X)\times\Jac_{\delta}(X)\times X$ be the
 map defined by 
 $\pi(D,L,x)=(\mu(D)\otimes K^{-1}\otimes L^{-1},L,x),$
 where $\mu:S^{2g-2+d}X\to\Jac_{2g-2+d}(X)$
 sends an effective divisor to the line bundle it represents.
 We have two Poincar\'e bundles $\LLL_1,\LLL_2\to
 \Jac_{\delta}(X)\times\Jac_{\delta}(X)\times X$ corresponding
 to each Jacobian.
 Let us define $\VVV=\pi^*\LLL_1\oplus\pi^*\LLL_2$. 

 On the other hand, we have a canonical bundle 
 $\LLL$ over $S^{2g-2+d}X\times X$ (the pullback of the Poincar\'e
 bundle through the projection $\mu\times\Id$) and a canonical section
 $\beta\in H^0(\LLL)$. Furthermore, we have an isomorphism
 $\pi_X^*K\otimes\pi^*\LLL_1\otimes\pi^*\LLL_2\cong p^*\LLL,$
 where $p:S^{2g-2+d}X\times\Jac_\delta(X)\times X\to S^{2g-2+d}X\times X$
 is the projection. Then
 we define $\Beta$ to be the canonical section $H^0(\pi_X^*K\otimes
 \pi^*\LLL_1\otimes\pi^*\LLL_2)$ induced by $\beta$.
 It is clear that $S_d$ is irreducible and that the family
 $(\VVV,\Beta)$ represents every isomorphism class of pairs
 of the form (\ref{eq:aquestsparells}).
 \end{pf}

 \begin{lemma}
 Any $\alpha$-polystable
 pair of the form (\ref{eq:aquestsparells}) can be deformed
 through a path of $\alpha$-stable pairs  
to a pair $(V,\beta)$ with $\Delta_{\beta}\neq 0$ and
 $\HH^0(\CCC_{V,\beta})=0$.
 \end{lemma}
 \begin{pf}
 By Lemmae \ref{totestable} and \ref{lemma:irrfamdelta}
 it suffices to prove that at least
 one pair of the form (\ref{eq:aquestsparells}) can be deformed
 to a pair $(V,\beta)$ with $\Delta_{\beta}\neq 0$ and
 $\HH^0(\CCC_{V,\beta})=0$.
 We will find two nonisomorphic
 line bundles $L_1,L_2$ of degree
 $\delta$ and nonzero sections
 $\beta\in H^0(K\otimes L_1\otimes L_2)$ and
 $\gamma\in H^0(K\otimes L_1^2)$.
 Then the family
 $\{(V_{\epsilon},\beta_{\epsilon})=
 (L_1\oplus L_2,\beta+\epsilon\gamma)\mid \epsilon\in\DD\}$
 is such desired deformation, since
 $(V_{\epsilon},\beta_{\epsilon})$
 is of the form (\ref{eq:aquestsparells})
 if and only if $\epsilon=0$.
 To find $L_1,L_2$, take two different points
 $\lambda,\lambda'\in\mu(S^{2g-2+d}X)$ and solve the
 following system of linear equations in $\Jac_{2g-2+d}(X)$:
(i) $K+L_1+L_2=\lambda$, (ii) $K+2L_1=\lambda'$.
 \end{pf}

 \subsection{Pairs of type C}
 \label{ss:typeC}

 Let $(V,\beta)$ be a pair with $\Delta=\Delta_{\beta}\neq 0$, and let
 $x=(x_1,\dots,x_k)$ be the vanishing locus of $\Delta$.
 Define the type of $(V,\beta)$ (with respect to the ordering of 
 the elements of $x$) to be
 $$\TTT=\TTT(V,\beta)=((r_1,z_1),\dots,(r_k,z_k))$$
 where $r_j$ (resp. $z_j$) is by definition the vanishing order
 of $\Delta$ (resp. $\beta$) at $x_j$.
 Note that by (\ref{estudilocal}) we always have $2z_j\leq r_j$.
 We also have
 \begin{equation}
 \sum_j r_j =\deg(K^2\otimes(\Lambda^2S)^2)=2(2g-2+d).
 \label{sumar}
 \end{equation}
 We define the generic type to be following list of $2(2g-2+d)$ pairs
 $$\TTT^{\gen}=((1,0),\dots,(1,0)).$$

 \begin{theorem}
 For any type $\TTT=((r_1,z_1),\dots,(r_k,z_k))$ there is a 
 manifold $S_{\TTT}$ with a covering $\bigcup U_i=S_{\TTT}$ and
 a family of pairs $(\VVV_i,\Beta_i)$ on $X$ parametrized by $U_i$
 for any $i$ so that

 {\rm (i)} if $s\in U_i\cap U_j$ then the pairs
 $(\VVV_i|_{\{s\}\times X},\Beta_i|_{\{s\}\times X})$
 and
 $(\VVV_j|_{\{s\}\times X},\Beta_j|_{\{s\}\times X})$
 are isomorphic;

 {\rm (ii)} if $\TTT=\TTT^{\gen}$ then $S_\TTT$ is irreducible;

 {\rm (iii)} any isomorphism class of
 pairs of type $\TTT$ which is $\alpha$-stable for some value of $\alpha$
 is represented by at least one
 $(V,\beta)=(\VVV_i|_{\{s\}\times X},\Beta_i|_{\{s\}\times X})$;

 {\rm (iv)} if at least one of the $r_j$ is odd, the
 dimension of $S_\TTT$ is $k+g-1+|x_R|/2$, where $|x_R|$ denotes
 the cardinal of $\{j\mid r_j\text{ is odd }\}$; this number
 is $\leq 7(g-1)+3d$, with equality if and only if $\TTT=\TTT^{\gen}$;

 {\rm (v)} if all $r_j$ are even then
 the dimension of $S_{\TTT}$ is strictly less than $7(g-1)+3d$.
 \end{theorem}

 \begin{pf}
 Let $k$ be a positive integer and define
 $B_\TTT=X^k\setminus\Delta_{\mult}$, where $\Delta_{\mult}$ 
 denotes the multidiagonal. 
 Let us fix a type $\TTT=((r_1,z_1),\dots,(r_k,z_k))$. 


 \noindent{\bf Case 1.} Suppose first that at least one $r_j$ is odd.

 Fix some $x=(x_1,\dots,x_k)\in B_\TTT$,
 let $x_R=\{x_j\in x\mid r_j\text{ is odd }\}$,
 and let $p:X'\to X$ be the $2:1$ covering ramified at $x_R$. 
 Take a pair $(V,\beta)$ with $\Delta_{\beta}^{-1}(0)=x$ and
 with type $\TTT$. Since $x_R\neq\emptyset$ there is no square root of
 $\Delta_{\beta}$.
 Let $(V',\beta')=p^*(V,\beta)$, and let
 $\Delta'=\Delta_{\beta'}$. By construction there is a square
 root of $\Delta'$, so there are line bundles
 $L_1,L_2\subset V'$ satisfying $\beta'\in H^0(p^*K\otimes L_1\otimes
 L_2)$.

 Let $\sigma:X'\to X'$ the Galois transformation. 
 There is a canonical lift $\sigma:V'\to V'$ which leaves
 $\beta'$ invariant and swaps $L_1$ and $L_2$
 (this follows from (i) in Lemma \ref{casD12}
 and the fact that if $\sigma^*L_1=L_1$ then $L_1$ would descend
 to a line bundle $L_0\to X$ satisfying 
 $\beta\in H^0(K\otimes L_0\otimes V)$, 
 hence $\Delta_{\beta}$ should have a square root).
 Denoting $L=L_1$ it then follows that
 $\beta'\in H^0(p^*K\otimes L\otimes\sigma^*L)$.

 By Lemma \ref{casD12} we have an exact sequence
 \begin{equation}
 0\to L\oplus\sigma^*L\to V'\to \cO_{T_x}\to 0,
 \label{extensioV}
 \end{equation}
 where $T_x$ is the following divisor on $X'$:
 $$T=\sum_{p(y)=x_j\in x\setminus x_R} (r_j/2-z_j) y+
 \sum_{p(y)=x_j\in x_R} (r_j-2z_j)y.$$
 In order for the sequence to be $\sigma$-equivariant
 one has to take the natural action of $\sigma$ on $\cO(T_x)$ obtained
 by identifying the sections of $\cO(T_x)$ with meromorphic functions
 on $X'$.
 So we have 
 $$2\deg L=\deg V'-\sum_j r_j+2\sum_j z_j =
 2d-\sum_j r_j+2\sum_j z_j.$$
 Consider the divisor 
 $$D_x=\sum_{p(y)=x_j\in x\setminus x_R} z_j y+
 \sum_{p(y)=x_j\in x_R} 2z_jy$$ in $X'$.
 Since the numbers
 $z_j$ describe the vanishing order of $\beta$, it follows that
 $p^*K\otimes L\otimes\sigma^*L\cong\cO(D_x)$.

 So begining from the pair $(V,\beta)$ we have constructed
 a line bundle $L$ on $X'$ of degree $d-\sum_j r_j/2+\sum_j z_j$
 satisfying $p^*K\otimes L\otimes\sigma^*L\cong\cO(D_x)$.

 Conversely, we can recover the isomorphism
 class of $(V,\beta)$ as follows. Let us denote 
 $$\Ext^1(\cO_{T_x},L\oplus\sigma^* L)_{\free}^{\sigma}\subset
 \Ext^1(\cO_{T_x},L\oplus\sigma^* L)^{\sigma}$$
 the set of elements giving locally free extensions.
 The $\sigma$-equivariant extensions of the type (\ref{extensioV})
 are classified by the elements of
 $\Ext^1(\cO_{T_x},L\oplus\sigma^* L)^{\sigma}$,
 and the isomorphism classes of
 $\sigma$-equivariant
 vector bundles $V'$ which are obtained through extensions
 of the form (\ref{extensioV}) are classified by
 $$\frac{\Ext^1(\cO_{T_x},L\oplus\sigma^* L)_{\free}^{\sigma}}
 {\Aut(L\oplus\sigma^* L)^{\sigma}\times\Aut(\cO_{T_x})}
 =\frac{\Ext^1(\cO_{T_x},L)_{\free}}{\Aut(L)\times\Aut(\cO_{T_x})},$$
 and the latter quotient has a unique element, represented by the 
 class of any extension
 \begin{equation}
 0\to L\to L\otimes\cO(T_x)\to \cO_{T_x}\otimes L\cong
 \cO_{T_x}\to 0.
 \label{extensioL}
 \end{equation}

 On the other hand, the section $\beta'$ (and hence $\beta$)
 is uniquely determined (up to multiplication by scalars)
 by the fact that it induces an isomorphism
 $p^*K\otimes L\otimes\sigma^*L\cong\cO(D_x)$.

 Let $d_\TTT=d-\sum_j r_j/2+\sum_j z_j$.
 We are now going to study the set
 $$\SSS_x:=\left\{L\in \Jac_{d_\TTT}(X')
 \mid p^* K\otimes L\otimes\sigma L\cong\cO(D_x)\right\}.$$
 The following lemma implies that 
 $\SSS_x$ is a torus of dimension $3(g-1)+d$.

 \begin{lemma}
 Let $\tau:\Jac_{d_\TTT}(X')\to \Jac_{2d_\TTT}(X')$
 be the map which sends $L$ to $L\otimes \sigma L$.
 The image of $\tau$ coincides with
 $\Jac_{2d_\TTT}(X')^{\sigma}$,
 and the fibres of $\tau$ can be identified to the torus
 $$\SSS_0=\left\{\Lambda\in \Jac_0(X')
 \mid \Lambda\otimes\sigma \Lambda\cong\cO\right\},$$
 which has dimension $g-1+|x_R|/2$.
\label{lemma:diminv}
 \end{lemma}
 \begin{pf}
 It is clear that the fibres of $\tau$ are the orbits of the
 action of $\SSS_0$ on $\Jac_{d_\TTT}(X')$ given by restricting the
 canonical action of $\Jac_0(X')$.

 For any $L\in\Jac_{2d_\TTT}(X')^{\sigma}$ there exists 
 a square root $L^{1/2}\in\Jac_{d_\TTT}(X')^{\sigma}$, which
 consequently satisfies $\tau(L^{1/2})=L$. It follows
 that $\Jac_{2d_\TTT}(X')^{\sigma}\subset \Im\tau$. The other
 inclusion is obvious, so we get 
 $\Jac_{2d_\TTT}(X')^{\sigma}=\Im\tau$.

 Let us now compute the dimension of $\SSS_0$.
 Recall that we have an identification 
 $\Jac_0(X')\cong  H^1(X';\RR)/H^1(X';\ZZ)$, which is
 $\sigma$-equivariant, since it comes
 from identifying $\Jac_0(X')$ with the set of gauge
 equivalence classes of flat $\U(1)$-connections on the
 trivial bundle on $X'$. Let $H^1(X';\RR)^{\pm}$ be the
 eigenspace of $\sigma$ corresponding to the
 eigenvalue $\pm 1$. We then have
 $T_{\uC}\SSS_0\cong  H^1(X';\RR)^-$ and
 $H^1(X;\RR)\cong H^1(X';\RR)^+$.
 So (using complex dimension everywhere)
 \begin{align*}
 \dim \SSS_0&=\dim H^1(X';\RR)^-
 =\dim H^1(X';\RR)-\dim H^1(X;\RR) \\
 &=g(X')-g=g-1+|x_R|/2.
 \end{align*}
 Here we have used Hurwitz' formula $\chi(X')=2\chi(X)-|x_R|$ to
 deduce $g(X')=1-\chi(X')/2=1-\chi(X)+|x_R|/2=-1+2g+|x_R|/2$
 \end{pf}

 Define $\pi:S_\TTT\to B_\TTT$ to be the fibration whose fibre
 over $x\in B_\TTT$ is $\SSS_x$. 
 As $x$ moves along $B_\TTT$ the divisors $x_R$ 
 sweep a divisor in $B_\TTT\times X$ which we denote by
 $\bx_R$. Let $p:X_\TTT'\to B_\TTT\times X$
 be the $2:1$ covering ramified along $\bx_R$, and let $\sigma:X_\TTT'\to
 X_\TTT'$ be the Galois transformation. 
 Let $D_\TTT$ (resp. $T_{\TTT}$)
 be the divisor in $X_\TTT'$ swept by the divisors $D_x$
 (resp. $T_x$) as $x$ moves along $B_\TTT$. 

 There is a universal bundle $\LLL\to S_\TTT\times_{B_\TTT}X_\TTT'$
 for which there is an isomorphism of bundles over
 $S_\TTT\times_{B_\TTT}X_\TTT'$
 \begin{equation}
 \cO(D_\TTT)\stackrel{\cong}{\longrightarrow}
 K\otimes\LLL\otimes\sigma^*\LLL
 \label{isoLLL}
 \end{equation}
 (we omit the pullbacks). 
 Now take a covering of $\SSS_{\TTT}$ by open sets $U_i$ admiting
 trivialisations
 $$\LLL\otimes \cO_{T_{\TTT}}|_{U_i\times X}
 \cong \cO_{T_{\TTT}}|_{U_i\times X}.$$

 Using each of these trivialisation we can obtain families of
 pairs on $U_i\times X$ beginning from the extensions (\ref{extensioL}).
 By the preceeding arguments it follows that any isomorphism class
 of pairs $(V,\beta)$ is represented in the family $(\VVV,\Beta)$.
 It is also clear that the dimension of $S_\TTT$ is $k+g-1+|x_R|/2$.
 Since $B_\TTT$ and the fibres of $S_\TTT\to B_\TTT$ are connected,
 we deduce that $S_\TTT$ is connected.

Finally, observe that $\TTT^{\gen}$ falls in the case we are now considering.
In this situation, $k=|x_R|=2(2g-2+d)$, so that
$$\dim \SSS_{\TTT^{\gen}}=7(g-1)+3d.$$
When $\TTT\neq\TTT^{\gen}$ then both $k$ and $|x_R|$ are
$\leq 2(2g-2+d)$, and one does not have equality in both
cases (i.e., $k=|x_R|=2(2g-2+d)$ only holds for $\TTT=\TTT^{\gen}$).

 \noindent{\bf Case 2.} Suppose that all the $r_j$ are even.
 Let $(V,\beta)$ be a pair of type $\TTT$ and vanishing locus $x\in B_\TTT$.
 There is a square root $\Delta^{1/2}$ of $\Delta_{\beta}$,
 so we have two line bundles $L_1,L_2$ on $X$, an exact sequence 
 $$0\to L_1\oplus L_2\to V\to\cO_{T_x}\to 0,$$
 with $T_x=\sum_j (r_j/2-z_j) x_j$ (see Lemma \ref{casD12}), and an isomorphism
 \begin{equation}
 \cO(D_x)\simeq K\otimes L_1\otimes L_2,
 \label{isoL12}
 \end{equation}
 where $D_x$ is the divisor $\sum_j z_j x_j$ in $X$.
 It follows that the isomorphism classes of pairs $(V,\beta)$
 with type $\TTT$ and vanishing locus $x$ are in 1---1
 correspondence with line bundles $L=L_1\in\Jac(X)$ (indeed, once 
 $L_1$ has been chosen, we set $L_2=\cO(D_x)\otimes K^{-1}\otimes
 L_1^{-1}$)
 and a choice of a class in
 \begin{equation}
 \frac{\Ext^1(\cO_{T_x},L_1\oplus L_2)_{\free}}
 {\Aut(L_1\oplus L_2)^{\sigma}\times\Aut(\cO_{T_x})}.
 \label{eq:extensiolliure}
 \end{equation}
 On the other hand, $\Ext^1(\cO_{T_x},L_1\oplus L_2)_{\free}$
 is the following open subset of $\Ext^1(\cO_{T_x},L_1\oplus L_2)$
 $$\Ext^1(\cO_{T_x},L_1\oplus L_2)_{\free}=
 \Ext^1(\cO_{T_x},L_1\oplus L_2)\setminus
 \bigcup_{f:\cO_{T_x}\exh\cO_{T}}
 f^*\Ext^1(\cO_T,L_1\oplus L_2),$$
 where the union is over all torsion sheaves $\cO_T$ admitting
 a surjection $f:\cO_{T_x}\exh\cO_{T}$ and
 $f^*:\Ext^1(\cO_T,L_1\oplus L_2)\to \Ext^1(\cO_{T_x},L_1\oplus L_2)$
 is the map induced by $f$.

 To compute $\Ext^1(\cO_{T_x},L_1\oplus L_2)$ we use Serre duality
 to obtain
 $\Ext^1(\cO_{T_x},L_1\oplus L_2)=\Hom(L_1\oplus L_2,\cO_{T_x}\otimes K),$
 and deduce from it that
 $$\dim \Ext^1(\cO_{T_x},L_1\oplus L_2)=2\sum r_j-2z_j.$$
 On the other hand, since we have 
 $\dim \Aut(\cO_T)=\sum r_j-2z_j$, we obtain:
 \begin{equation}
 \dim \frac{\Ext^1(\cO_{T_x},L_1\oplus L_2)_{\free}}
 {\Aut(L_1\oplus L_2)^{\sigma}\times\Aut(\cO_{T_x})}=
 \left(\sum r_j-2z_j\right)-2\leq\sum r_j-2=2(2g-2+d)-2.
 \label{eq:dimExt}
 \end{equation}

 Now, if the pair $(V,\beta)$ is $\alpha$-polystable for some value of $\alpha$
 then we necessarily have $\deg L_1\leq d/2$ and $\deg L_2\leq d/2$.
 (\ref{isoL12}) implies that 
 \begin{align*}
 \deg L_1 &= \deg(D_x)-\deg K-\deg L_2 \\
 &= \sum z_j+2-2g-\deg L_2,
 \end{align*}
 so semistability implies that
 $$\deg L_1\in [m_\TTT,M_\TTT]:=\left[\sum z_j+2-2g-d/2,d/2\right].$$
 Consequently we can define
 $$S_\TTT:=\bigsqcup_{m_\TTT\leq d\leq M_\TTT}\EEE^d(X),$$
 where $\EEE^d\to\Jac_d(X)\times B_{\TTT}$ is the bundle whose fibre over 
 $(L,x)$ is
 $$\frac{\Ext^1(\cO_{T_x},L\oplus \cO(D_x)\otimes K^{-1}\otimes L^{-1})}
 {\Aut(L\oplus \cO(D_x)\otimes K^{-1}\otimes L^{-1})\times\Aut(\cO_{T_x})}.$$
 By (\ref{eq:dimExt}) this has dimension
 $\dim\EEE^d\leq 5(g-1)+2d+k-1.$
 Finally, since all zeros of $\Delta_{\beta}$ have even order (hence 
 $\geq 2$) we have $k\leq 2(g-1)+d=\deg(K^2\otimes(\Lambda^2V)^2)/2$ so
 $\dim\EEE^d\leq 7(g-1)+3d-1.$
 The definitions of the local families $\VVV_i$ and $\Beta_i$
 can be given exactly as in Case 1.
 \end{pf}

\subsection{Existence of stable objects}
\label{ss:existencia}

Here we prove that for any $-2(g-1)<d\leq 0$ and $\alpha\geq 0$
there is a $\alpha$-stable pair $(V,\beta)$ with $\deg V=0$. More concretely,
we show that there exists a $\alpha$-stable pair of degree $d$
and generic type, i.e., such that $\Delta_{\beta}$
has simple zeroes (and consequently $\beta$ never vanishes).
By the results in Subsection \ref{ss:typeC}, the isomorphism
classes of pairs of 
degree $d$ of generic type are in 1---1 correspondence with
choices of:
\begin{itemize}
\item a $(2:1)$ covering $p:X'\to X$, ramified along
$2(2g-2+d)$ different points of $X$;
\item a line bundle $L$ on $X'$ of degree $-2(g-1)$ such that 
\begin{equation}
p^*K\otimes L\otimes \sigma^*L\cong\uC,
\label{eq:trivL}
\end{equation}
where
$\sigma:X'\to X'$ is the unique nontrivial
automorphism of $p$.
\end{itemize}
The pair $(V,\beta)$ corresponding to one such choice is related
to this data as follows: there is an exact sequence
\begin{equation}
0\to L\oplus\sigma^*L\to p^*V\to \cO_R\to 0,
\label{eq:relLV}
\end{equation}
where $\cO_R$ is the structure sheaf of the ramification
locus, and $p^*\beta\in H^0(p^*K\otimes L\otimes\sigma^*L)$
is a nonvanishing section giving rise to the isomorphism
(\ref{eq:trivL}).

Define, for every $k\in\ZZ$, 
$\Jac_k^+:=\{L\in\Jac_k(X')\mid\sigma^*L\cong L\}$
and $\Jac_k^-:=\{L\in\Jac_k(X')\mid\sigma^*L\cong L^{-1}\}$.
It is easy to check that for any pair $p,q\in\ZZ$ the map
given by tensor product
$\Jac_p^+\times\Jac_q^-\to\Jac_{p+q}(X')$
is a covering map.

Pick some integer $\delta\geq-2(2g-2+d)$ and define
$$W_{\delta}:=\{\LLL\in\Jac_{-2(g-1)-\delta}(X')\mid h^0(\LLL)>0\}.$$
This is a complex submanifold of $\Jac_{-2(g-1)-\delta}(X')$,
and we have 
\begin{equation}
\dim W_{\delta}\leq -2(g-1)-\delta.
\label{eq:boundW}
\end{equation}
Take a square root $K^{1/2}\in\Jac(X)$ of $K$.
Let also $Q$ be the quotient $\Jac_{-2(g-1)-\delta}/\Jac_{-\delta}^+$,
and consider the diagram
$$\xymatrix{ & W_{\delta}\ar[d]_{c} \\
\Jac_0^-\ar[r]^{f} & Q,}$$
where $\alpha$ maps $\LLL$ to its class $[\LLL]$ in $Q$ and
$f$ maps $L_0\in\Jac_0^-$ to $[p^*K^{1/2}\otimes L_0]$.
The map $f$ is a covering, hence $\dim Q=\dim\Jac_0^-$, and
this dimension is equal to 
$$\dim T_{\uC}\Jac_0^-=\dim H^1(X';\RR)^-=3(g-1)+d.$$
(see the proof of Lemma \ref{lemma:diminv}).
It is then straightforward to check that for any 
$\delta\geq-2(2g-2+d)$ we have $\dim W_{\delta}<\dim Q$.
So, if we set 
$$W:=\bigcup_{\delta \geq-2(2g-2+d)}W_{\delta},$$
which is a finite union because of (\ref{eq:boundW}), then 
$c(W)\neq Q$. Hence there is some $L_0\in\Jac_0^-$ such that
$f(L_0)\notin c(W)$. Let $L:=p^*K^{1/2}\otimes L_0$, and let
$(V,\beta)$ be the pair constructed using the ideas in 
Subsection \ref{ss:typeC}.

We now prove that $(V,\beta)$ is $\alpha$-stable.
First of all, by construction $\Delta_{\beta}$
does not admit a square root. Hence, $(V,\beta)$ could only 
be unstable if there were a line subbundle $\Lambda\subset V$
such that $\deg \Lambda\geq c\geq 0$.
In this case, by (\ref{eq:relLV}), we would have a diagram
$$\xymatrix{0 \ar[r] & L\oplus\sigma^*L \ar[r]
& p^*V \ar[r] & \cO_R \ar[r] & 0 \\
0 \ar[r] & M \ar[u]^{(j,j')}\ar[r] & p^*\Lambda \ar[u]\ar[r]
& T \ar[u]^{i}\ar[r] & 0,}$$
in which $T$ is a torsion sheaf which injects into $\cO_R$ by $i$.
Consequently, we have
$$\delta:=\deg M\geq \deg p^*\Lambda-|\cO_R|=2c-2(2g-2+d)\geq-2(2g-2+d).$$
On the other hand, $M\in\Jac_{\delta}^+$ and the map $(j,j')$
is $\sigma$ invariant, hence $j'=\sigma^*j$. And, since
$(j,j')$ is an inclusion, it follows that $j:M\to L$
is a nonzero holomorphic map. Hence $\LLL:=M^*\otimes L\in W_{\delta}$.
But this implies that $f(L_0)=[L]\in c(W_{\delta})$, contradicting
our assumption. So there does not exist any subbundle
$\Lambda\subset V$ with $\deg\Lambda\geq 0$, thus $(V,\beta)$ is
$\alpha$-stable.

\noindent {\bf Acknowledgements}.
We wish  to thank Nigel Hitchin for very useful discussions. 
Much of the research  in this paper was done while the
authors were visiting  the Mathematical Institute of Oxford
(which we warmly thank for its hospitality) during
Michaelmas term of 2001. We thank
the EPSRC and EDGE for supporting the visists of the first and  second
authors respectively.
The  authors are members of VBAC (Vector Bundles on Algebraic Curves),
which is partially supported by EAGER (EC FP5 Contract no. HPRN-CT-2000-00099)
and by EDGE (EC FP5 Contract no. HPRN-CT-2000-00101). 


 \end{document}